%% file: Hecke.tex
\documentclass[11pt, a4paper]{article}

\usepackage{amssymb}    
\usepackage{latexsym} 
\usepackage{graphicx}
\usepackage{html}       

\newcommand{\C}{\ensuremath{\mathbb{C}}}
\newcommand{\Q}{\ensuremath{\mathbb{Q}}}
\newcommand{\R}{\ensuremath{\mathbb{R}}}
\newcommand{\Z}{\ensuremath{\mathbb{Z}}}
\newcommand{\N}{\ensuremath{\mathbb{N}}}

\newcommand{\FE}{\ensuremath{\mathrm{FE}}}  
\newcommand{\LL}{\ensuremath{\mathcal{L}}}

\newcommand{\RR}{\ensuremath{\mathcal{R}}}

\newcommand{\HH}{\ensuremath{\mathfrak{H}}} 
\newcommand{\HHstar}{\ensuremath{\mathfrak{H}^\ast}}

\newcommand{\dd}{\ensuremath{\mathrm{d}}}

\newcommand{\sign}[1]{\ensuremath{\mathrm{sign}\left(#1\right)}}

\newcommand{\Gmod}{\ensuremath{{\Gamma(1)}}}
\newcommand{\Gnull}[1]{\ensuremath{{\Gamma_0\left(#1\right)}}}

\newcommand{\re}[1]{\ensuremath{{\mathrm{Re}\!\left( #1 \right)}}}
\newcommand{\im}[1]{\ensuremath{{\mathrm{Im}\!\left( #1 \right)}}}
\newcommand{\SL}[1]{\ensuremath{{\mathrm{SL}_2\!\left( #1 \right)}}}

\newcommand{\GL}[1]{\ensuremath{{\mathrm{GL}_2\!\left( #1 \right)}}}

\newcommand{\Matrix}[4]{{\textstyle \left( {#1\atop #3} \: {#2 \atop  #4} \right)}}

\newtheorem{theorem}{Theorem}[section]
\newtheorem{lem}[theorem]{Lemma}
\newtheorem{prop}[theorem]{Proposition}

\newtheorem{defn}[theorem]{Definition}

\newenvironment{proof*}[1][]{\smallskip\par\noindent\emph{Proof}#1{\emph. }
 \ignorespaces}{\nopagebreak\hfill$\Box$\smallskip\vspace{1ex}\par\ignorespaces}

\begin{document}

\title{Hecke operators on period functions for $\Gamma_0(n)$}
\author{T. M\"uhlenbruch\\
        \small    Institut f\"ur Theoretische Physik\\
        \small    TU Clausthal, Germany\\
        \small    \texttt{tobias.muehlenbruch@tu-clausthal.de}}
\date{May 10, 2005}
\maketitle

\begin{abstract}
Matrix representations of Hecke operators on classical holomorphical cusp forms and the corresponding period polynomials are well known. 
In this article we derive representations of Hecke operators for vector valued period functions for the congruence subgroups $\Gnull{n}$. 
For this we use an integral transform from the space of vector valued cusp forms to the space of vector valued period functions.\\[1ex]
\end{abstract}

\input{introduction}

\input{farey}

\input{cusp}

\input{period}

\input{results}

\textbf{Acknowledgments.} 
This work has been supported by the Deutsche Forschungsgemeinschaft through the DFG Forschergruppe ``Zetafunktionen und lokalsymmetrische R\"aume''.
The author thanks Dieter Mayer and Joachim Hilgert for their discussions and comments.

\def\cprime{$'$}

\bigskip
\noindent
{\small
Tobias M\"uhlenbruch\\
Institut f\"ur Theoretische Physik\\
Abteilung Statistische Physik und Nichtlineare Dynamik\\
Technische Universit\"at Clausthal\\
Arnold-Sommerfeld-Stra{\ss}e 6\\
38678 Clausthal-Zellerfeld\\
Germany\\[1ex]
email: tobias.muehlenbruch@tu-clausthal.de
}

\end{document}

%% file: introduction.tex
\section{Introduction}
\label{A}
In \cite{muehlenbruch:13} we discussed the Hecke operators for period functions for the full modular group.
In particular we showed that the formal sum of integral matrices $\sum_{{a>c \geq 0 \atop d>b \geq 0} \atop ad-bc=m } \Matrix{a}{b}{c}{d}$ is a representation of the $m^\mathrm{th}$ Hecke operator on period functions for the full modular group.
The aim of the present paper is to give an explicit representation of the $m^\mathrm{th}$ Hecke operator on period functions for subgroups $\Gnull{n}$.
The method used in this paper holds for all discrete subgroups of the full modular group with finite index.

\smallskip

To state our main result and to sketch the content of each section we have to fix some notations used throughout the text. 
For an integer $n$ let $\mathrm{Mat}_n(2,\Z)$ (respectively $\mathrm{Mat}_\ast(2,\Z)$) be the set of $2\times 2$ matrices with integer entries and determinant $n$ (respectively nonzero determinant). 
Let $\RR_n:= \Z[\mathrm{Mat}_n(2,\Z)]$ (respectively $\RR:= \Z[\mathrm{Mat}_\ast(2,\Z)]$) be the set of finite linear combinations with coefficients in $\Z$ of elements of $\mathrm{Mat}_n(2,\Z)$ (respectively $\mathrm{Mat}_\ast(2,\Z)$).
Similarly, we denote by $\mathrm{Mat}_n^+(2,\Z)$ (respectively $\mathrm{Mat}_\ast^+(2,\Z)$) the set of $2\times 2$ matrices with integer nonnegative entries and determinant $n$ (respectively nonzero determinant). 
Denote furthermore by $\RR_n^+:= \Z[\mathrm{Mat}_n^+(2,\Z)]$ (respectively $\RR^+:= \Z[\mathrm{Mat}_\ast^+(2,\Z)]$) the set of finite linear combinations with coefficients in $\Z$ of elements of $\mathrm{Mat}_n^+(2,\Z)$ (respectively $\mathrm{Mat}_\ast^+(2,\Z)$).
Note that $\RR^{(+)} = \bigcup_{n=1}^\infty \RR_n^{(+)}$ and $\RR_n^{(+)} \cdot \RR_m^{(+)} \subset \RR_{nm}^{(+)}$.
By definition we have $\SL{\Z} = \mathrm{Mat}_1(2,\Z)$.
The following four elements of $\SL{\Z}$ will play a special role in our paper:
\[
I:=       \Matrix{1}{0}{0}{1}, \quad
T:=       \Matrix{1}{1}{0}{1}, \quad
S:=       \Matrix{0}{-1}{1}{0} \quad \mbox{and} \quad
T^\prime:=\Matrix{1}{0}{1}{1}.
\]
For each $n \in \N$ the Hecke congruence subgroup $\Gnull{n}$ is given by 
\[
\Gnull{n} :=\left\{
\Matrix{a}{b}{c}{d} \in \SL{\Z};\, c \equiv 0 \bmod n
\right\}.
\]
The full modular group is denoted by $\Gmod:=\Gnull{1}=\SL{\Z}$. 
Let $\mu=\mu_n$ denote the index of $\Gnull{n}$ in $\Gmod$ and let $\alpha_1,\ldots,\alpha_\mu$ denote representatives of the right cosets in $\Gnull{n}\backslash \Gmod$.

We also need the set of upper triangular matrices
\begin{equation}
\label{A1.2}
X_m =
\left\{ \Matrix{a}{b}{0}{d} \in \mathrm{Mat}_m(2,\Z); \, d > b \geq 0 \right\}.
\end{equation}

Let $\bar{\C}$ denote the one-point compactification of $\C$.
The map
\[
\mathrm{Mat}_\ast(2,\Z) \times \bar{\C} \to \bar{\C};\quad
\left(\Matrix{a}{b}{c}{d},z\right) \mapsto \Matrix{a}{b}{c}{d}z :=\frac{az+b}{cz+d}.
\]
gives an action of the matrices on $\bar{\C}$.
It induces the familiar \emph{slash action} $|_s$ on functions $f$ on $\HH:=\{z\in \C;\, \im{z}>0\}$ resp.\ $(0,\infty)$ formally defined by
\begin{equation}
\label{A1.1}
f\big|_s\Matrix{a}{b}{c}{d}(z) =
(ad-bc)^s\, (cz+d)^{-2s} \,f\left(\frac{az+b}{cz+d}\right)
\end{equation}
for complex numbers $s$ and certain classes of matrices $\Matrix{a}{b}{c}{d}$.
The slash action is well defined for
\begin{description}
\item[(a)]
$s \in \Z$, $\Matrix{a}{b}{c}{d} \in \mathrm{Mat}_m(2,\Z)$, $m\in \N$, $z \in \HH$ and
\item[(b)]
$s \in \C$, $\Matrix{a}{b}{c}{d} \in \mathrm{Mat}_m^+(2,\Z)$, $m\in \N$ and $z \in (0,\infty)$
\end{description}
as the discussion in~\cite{hilgert:1} shows.
A simple calculation gives $\big(f\big|_s\alpha\big) \big|_s \gamma = f\big|_s(\alpha\gamma)$ for all matrices $\alpha$, $\gamma \in \mathrm{Mat}_\ast(2,\Z)$ in case (a) resp.\ $\mathrm{Mat}_\ast^+(2,\Z)$ in case (b).
We extend the slash action linearly to formal sums of matrices.

\smallskip

Recall Maass cusp forms: 
A \emph{cusp form} $u$ for the congruence subgroup $\Gnull{n}$ is a real-analytic function $u:\HH \to \C$ satisfying:
\begin{enumerate}
\item
$u(gz) =u(z)$ for all $g \in \Gnull{n}$,
\item
$\Delta u= s(1-s)u$ for some $s \in \C$ where $\Delta = -y^2(\partial_x^2+\partial_y^2)$ is the hyperbolic Laplace operator.
We call the parameter $s$ the \emph{spectral parameter} of $u$.
\item
\label{A1.3}
$u$ is of rapid decay in all cusps:
if $p \in \Q \cup \{\infty\}$ is a cuspidal point for $\Gnull{n}$ and $g\in\Gmod$ is such that $gp=\infty$ then $u(gz) = \mathrm{O}\left(\im{z}^C\right)$ as $\im{z} \to \infty$ for all $C \in \R$.
\end{enumerate} 
We denote the space of cusp forms for $\Gnull{n}$ with spectral value $s$ by $S(n,s)$.

The last condition above actually states two different conditions: vanishing in all cusps and an explicit growth condition.
Both are equivalent since it is shown in \cite{iwaniec:1} that vanishing in the cusp $p=g^{-1} \infty$ implies the stronger growth condition $u(gz) = \mathrm{O}\left(e^{-\im{z}}\right)$ as $\im{z} \to \infty$ which again implies vanishing at the cusp $p=g^{-1}\infty$.

\smallskip

A function $f:(0,\infty) \to \C$ is called \emph{holomorphic} if it is locally the restriction of a holomorphic function.

The vector valued period functions for $\Gnull{n}$ are defined as follows:
A \emph{period function} for $\Gnull{n}$ is a function $\vec{\psi}:(0,\infty) \to \C^\mu$ with $\vec{\psi}=(\psi_i)_{i\in\{1,\ldots,\mu\}}$ such that
\begin{enumerate}
\item $\psi_i$ is holomorphic on $(0,\infty)$ for all $i \in \{1\ldots,\mu\}$.
\item $\vec{\psi}(z) = \rho(T^{-1})\, \vec{\psi}(z+1) + (z+1)^{-2s}\,\rho({T^\prime}^{-1})\, \vec{\psi}\left(\frac{z}{z+1}\right)$.
We call the parameter $s\in \C$ the \emph{spectral parameter} of $\vec{\psi}$.
The matrix representation $\rho:\Gmod\to \C^{\mu\times \mu}$ is induced by the trivial representation of $\Gnull{n}$ as introduced in \S\ref{D1}.
\item 
For each $i=1,\ldots,\mu$ $\psi_i$ satisfies the growth condition 
\[
\psi_i(z) = \left\{\begin{array}{ll}
\mathrm{O} \left( z^{\max\{0,-2\re{s}\}} \right) \quad& \mbox{as } z \downarrow 0 \mbox{ and}\\
\mathrm{O} \left( z^{\min\{0,-2\re{s}\}} \right) \quad& \mbox{as } z\to \infty.
\end{array} \right.
\]
\end{enumerate} 
Following \cite{lewis:2} we denote the space of period functions for $\Gnull{n}$ with spectral value $s$ by $\FE(n,s)$.
We call a function $\vec{\psi}$ a \emph{period like function} if $\vec{\psi}$ satisfies only the conditions $1$ and $2$. 
The space of period like functions for $\Gnull{n}$ with spectral value $s$ is denoted by $\FE^\ast(n,s)$

\smallskip

In \S\ref{D} we construct a bijective map $S(n,s) \to \FE(n,s)$:
For this we introduce the space of vector valued cusp forms $S_\mathrm{ind}(n,s)$ and define an operator $P:S_\mathrm{ind}(n,s) \to \FE(n,s)$ in the form of an explicit integral transform.
The operator $P$ was studied in \cite{lewis:2} for cusp forms of the full modular group and in~\cite{muehlenbruch:12} for cusp forms of $\Gmod$ with arbitrary real weight.
Martin discusses a similar integral operator for modular cusp forms of weight 1 in \cite{martin:1}.

We derive our main result through a map $S(n,s) \to \FE(n,s)$ by which the \linebreak $m^\mathrm{th}$ Hecke operator $H_{n,m}$ on $S(n,s)$ induces an operator $\tilde{H}_{n,m}$ on $\FE(n,s)$.
For $n=1$ and $m$ prime these operators coincide with the operators $\tilde{T}(m)$ in \cite{hilgert:1} which are derived there from the transfer operators only. 
For this special case the operators were also determined in \cite{muehlenbruch:13} using Eichler integrals to map Hecke operators on cusp forms to Hecke operators on period functions.
For $n >1$ the operators $\tilde{T}_{n,m}$ in~\cite{hilgert:1} and the operators $\tilde{H}_{n,m}$ are not the same in general. 
The exact relation between these operators will be discussed in \cite{mayer:8}.

%% file: farey.tex
\section{On Farey sequences}
\label{C}
In this section we discuss some properties of rational numbers related to Farey sequences.
Our discussion is closely related to the modified continued fractions introduced in \cite{hilgert:1}.

\subsection{Farey sequences a la Hurwitz}
\label{C1}
Let us recall the theory of Farey-sequences. 
Most of the properties mentioned can be found in \cite{hurwitz:1}. 
We adhere to the convention to denote infinity in rational form as $\infty = \frac{1}{0}$ and $-\infty = \frac{-1}{0}$ and to denote rationals $\frac{p}{q}$ with coprime $p \in \Z$ and $q \in \N$.

\begin{defn}
\label{C1.1}
For $n \in \N$ the \emph{Farey-sequence} $F_n$ of \emph{level $n$} is the sequence 
\[
F_n := \left( \frac{u}{v};\, u,v \in \Z, |u| \leq n, 0 \leq v \leq n \right).
\]
ordered by the standard order $<$ of $\R$.
We define $F_0$ as
\[
F_0 := \left( \frac{-1}{0}, \frac{0}{1}, \frac{1}{0} \right).
\]
The level function $\mathrm{lev}: \Q \rightarrow \Z$ is defined by
\[
\mathrm{lev}\left(\frac{a}{b}\right)=\left\{ \begin{array}{ll} 
0 & \mbox{if } \frac{a}{b} \in \{ \frac{-1}{0}, \frac{0}{1}, \frac{1}{0} \} \mbox{ and}\\
\max\{|a|,|b|\}\qquad & \mbox{otherwise}.
\end{array}\right.
\]
\end{defn}

Let $\frac{a}{c}$ and $\frac{b}{d}$ be two neighbors in the Farey-sequence $F_n$. 
Then the square matrix $\Matrix{a}{b}{c}{d}$ satisfies
$\det \Matrix{a}{b}{c}{d} = \pm 1$ (see Satz 1 in \cite{hurwitz:1}).

\begin{lem}
\label{C1.2}
Let $\frac{a}{c}$ and $\frac{b}{d}$ be two neighbors of the Farey-sequence $F_n$.
Then 
\[
\det \Matrix{a}{b}{c}{d} = -1 
\quad \Longleftrightarrow \quad \frac{a}{c} < \frac{b}{d}.
\]
\end{lem}

\begin{proof*}
For $c=0$ (resp.\ $d=0$) we have $\frac{-1}{0} < \frac{-n}{1}$ (resp.\ $\frac{n}{1}<\frac{1}{0}$) and $\det \Matrix{-1}{-n}{0}{1} = -1$ (resp.\ $\det \Matrix{n}{1}{1}{0} = -1$).

Assume that $c,d>0$.
Since $\frac{a}{c} < \frac{b}{d}$ is equivalent to $\frac{ad}{cd} < \frac{bc}{cd}$ and hence to $ad-bc<0$, the statement of the lemma follows from Satz~1 in~\cite{hurwitz:1}.
\end{proof*}

\textit{Remark.} 
Our applications of the Farey sequences deal mostly with the case $\det \Matrix{a}{b}{c}{d} = -1$. 
However, we prefer matrices in $\Gmod$. 
For this we replace the matrix $\Matrix{a}{b}{c}{d}$ by $\Matrix{-a}{b}{-c}{d}$. 
This obviously does not change the rational numbers $\frac{a}{c}$ and $\frac{b}{d}$.

\smallskip 

We need also the following result in \cite{hurwitz:1}:
\begin{lem}
\label{C1.3}
For $\frac{a}{c}$ and $\frac{b}{d}$ with $a,b,c,d\in \Z$, $c,d \geq 0$ and $ad-bc = \pm 1$ define $n:=\max \left\{ \mathrm{lev}\left(\frac{a}{c}\right),\, \mathrm{lev}\left(\frac{b}{d}\right) \,\right\}$.
Then $\frac{a}{c}$ and $\frac{b}{d}$ are neighbors in the Farey-sequence $F_n$.
\end{lem}

\subsection{Left neighbor sequences}
\label{C2}
We define the \emph{left neighbor map} $\mathrm{LN}: \Q \cup \{+\infty\} \to  \Q \cup \{-\infty\}$ such that $\mathrm{LN}(q)$ is the left neighbor of $q$ in the Farey sequence $F_{\mathrm{lev}(q)}$, that is
\begin{equation}
\label{C2.1}
\mathrm{LN}(q) = 
\max\{ r \in F_{\mathrm{lev}(q)}; \, r < q \}.
\end{equation}
For $q=\frac{b}{d} \in \Q \cup \{+\infty\}$ put $\frac{a}{c}= \mathrm{LN}\big(\frac{b}{d}\big) \in \Q \cup \{-\infty\}$.
Then by construction of the map $\mathrm{LN}$ we have $\mathrm{lev}\big(\frac{a}{c}\big) \leq \mathrm{lev}\big(\frac{b}{d}\big)$ and by Lemma~\ref{C1.2} $\det \Matrix{a}{b}{c}{d} = -1$.

\begin{lem}
\label{C2.2}
For $q \in \Q\cup \{+\infty\}$ and $\mathrm{lev}(q) >0$ we have $\mathrm{lev}\big(\mathrm{LN}(q)\big) < \mathrm{lev}(q)$.
\end{lem}

\begin{proof*}
If $\mathrm{lev}(q)=1$ the set of rational numbers of level $1$ is $\{\pm1\}$.
The statement of the lemma then follows since $\mathrm{LN}(1)=0$ and $\mathrm{LN}(-1)=-\infty$ are of level $0$.

If $N=\mathrm{lev}(q)>1$ write $q=\frac{b}{d}$ with $\gcd(b,d)=1$, $d \geq 0$ and $\mathrm{LN}(q)=\frac{a}{b}$ with $\gcd(a,b)=1$, $b \geq 0$. Lemma~\ref{C1.2} implies that $ad-bc=-1$. 
Assume that $\mathrm{lev}\big(\mathrm{LN}(q)\big) = \mathrm{lev}(q)$. 
Then $\max\{|a|,c\}=\max\{|b|,d\}=N$ and $|a| \not= c$, $|b| \not= d$.
There are four cases to consider:
\begin{description}
\item[if $|a|=|b|=N$] 
  we find $-1=ad-bc= \sign{a}Nd-\sign{b}Nc$. 
  Hence, $\sign{b}c-\sign{a}d=\frac{1}{N}$ which contradicts $c,d\in\Z$,
\item[if $c=d=N$] 
  we find $-1=ad-bc= aN-bN$. Hence, $b-a=\frac{1}{N}$ which contradicts $a,b\in\Z$.
\item[if $|a|=d=N$]
  we find $-1=ad-bc= \sign{a}N^2-bc$. Hence, $bc=1+\sign{a}N^2$ and $|b|,c \leq N-1$.
  The second estimate implies $|bc| <N^2-1$ which is a contradiction.
\item[if $|b|=c=N$] 
  we find $-1=ad-bc= ad-\sign{b}N^2$. Hence, $ad=-1+\sign{b}N^2$ and $|a|,d \leq N-1$.
  The second estimate implies $|ad| <N^2-1$ which is a contradiction.
\end{description} 
Hence the assumption $\mathrm{lev}\big(\mathrm{LN}(q)\big) = \mathrm{lev}(q)$ was wrong and $\mathrm{lev}\big(\mathrm{LN}(q)\big) < \mathrm{lev}(q)$ must hold.
\end{proof*}

\begin{defn}
\label{C2.3}
Let be $q \in \Q \cup \{+\infty\}$ and $L=L_q\in \N$ such that
\[
\mathrm{LN}^{L}(q) = -\infty
\qquad \mbox{and} \qquad
\mathrm{LN}^{l}(q) > -\infty
\quad \mbox{for all } l=1,\ldots,L-1.
\]
The \emph{left neighbor sequence} $\mathrm{LNS}(q)$ of $q$ is the finite sequence
\[
\mathrm{LNS}(q)= \big( \mathrm{LN}^{L}(q), \mathrm{LN}^{L-1}(q), \ldots, \mathrm{LN}^{1}(q), q \big),
\]
where we use the notation $\mathrm{LN}^{l}(q) := \underbrace{\mathrm{LN}\big(\mathrm{LN}(\cdots \mathrm{LN}(q)) \cdots \big)}_{l\,\mathrm{times}}$. 
\end{defn}

\textit{Remark.}
The number $L$ in Definition~\ref{C2.3} is unique. 

\begin{lem}
\label{C2.4}
For $q\in \Q$ consider the left neighbor sequence 
\[
\mathrm{LNS}(q)=(y_0, \ldots, y_L)
\]
with $y_0=-\infty$ and $y_L=q$.
The sequence $(y_L,\ldots,y_0)$ is a \emph{partition} of $q$ which is \emph{minimal} for $0 \leq q < 1$ in the sense of Definition~2.3 in~\cite{hilgert:1}.
\end{lem}

\begin{proof*}
We assume the elements of $\mathrm{LNS}(q)$ to be given as $y_l = \frac{a_l}{b_l}$ with $\gcd(a_l,b_l)=1$ and $b_l\geq 0$. 
By construction, $y_{l-1}< y_l$ and both numbers are neighbors in the Farey sequence $F_{\mathrm{lev}(y_l)}$.
Lemma~\ref{C1.2} implies $\det\Matrix{a_{l-1}}{a_l}{b_{l-1}}{b_l} = -1$.
Hence, $\det\Matrix{b_l}{-a_l}{b_{l-1}}{-a_{l-1}} = 1$ and the sequence $(y_L,\ldots,y_0)$ is a partition of $q$ in the sense of Definition~2.3 in~\cite{hilgert:1}. 

We consider the case $0 \leq q < 1$. 
We have to show that the partition is minimal in the sense of Definition~2.3 in~\cite{hilgert:1} and hence the denominators $b_l$ of $y_l$ have to satisfy
\[
0=b_0 < b_1 < \ldots < b_{L-1} < b_L.
\]
The construction of the left neighbour map and the fact $q >0$ implies that $y_1=0$ and $b_1=1$.
Since $0=y_1< y_l <1$ for all $l =2,\ldots,L$ forces the denominator $b_l$ of $y_l$ to be larger than $1$ for all $l=2,\ldots,L$ we have to check the inequalities above only for the indices $l\geq 2$. 
Obviously for $0<\frac{a}{b}<1$ with $\gcd(a,b)=1$ and $a,b\geq 0$ one has $b>a$ and hence $\mathrm{lev}\big(\frac{a}{b}\big) = b$.

Consider $y_l$, $l=2,\ldots,L$, with denominator $b_l= \mathrm{lev}\big(\frac{a_l}{b_l}\big) >0$. 
Lemma~\ref{C2.2} implies that $y_{l-1}=\mathrm{LN}(y_l)$ satisfies $b_{l-1}= \mathrm{lev}(y_{l-1})< b_l$. 
The partition $(y_L,\ldots,y_0)$ is indeed minimal. 
\end{proof*}

\begin{lem}
\label{C2.5}
For $0<q<1$ rational the two sequences $\mathrm{LNS}(q)=(y_0, \ldots, y_L)$ with $y_0=-\infty$ and $y_L=q$ and the sequence defined by the modified continued fraction expansion $(x_0, \ldots, x_{L^\prime})$ with $x_0=q$ and $x_{L^\prime}=-\infty$ given in \cite{hilgert:1} coincide.
Indeed $L=L^\prime$ and $y_l=x_{L-l}$ for all $l=0,\ldots,L$.
\end{lem}

\begin{proof*}
For $0\leq q<1$ rational the left neighbor sequence $\mathrm{LNS}(q)$ is a minimal partition by Lemma~\ref{C2.4}. 
According to \cite{hilgert:1} this partition is unique and hence Lemma~\ref{C2.5} holds.
\end{proof*}

\begin{defn}
\label{C2.6}
To $q \in [0,1)$ rational and $\mathrm{LNS}(q)=\big(\frac{a_0}{b_0}, \ldots,\frac{a_L}{b_L} \big)$ with 
\linebreak
$\gcd(a_l,b_l)=1$ and $b_l \geq 0$, $l=0,\ldots,L$, we attach the element $M(q)=\sum_{l=1}^L m_l \in \RR_1$:
\begin{equation}
\label{C2.7}
M(q) = \Matrix{-a_0}{a_1}{-b_0}{b_1}^{-1} + \ldots + 
       \Matrix{-a_{l-1}}{a_{l}}{-b_{l-1}}{b_{l}}^{-1} + \ldots +
       \Matrix{-a_{L-1}}{a_{L}}{-b_{L-1}}{b_{L}}^{-1}.
\end{equation} 
\end{defn}

\textit{Remark.} 
Obviously the number $L$ in Definition~\ref{C2.6} depends on $q$.

\smallskip

In the following we need some properties of the matrices in $M(q)A$ for $A \in \mathrm{Mat}_\ast(2,\Z)$:

\begin{lem}
\label{C2.8}
For $0\leq q < 1$ rational and $M(q)=\sum_{l=1}^L \Matrix{\ast}{\ast}{c_l}{d_l}$ one has $c_l \zeta +d_l >0$ for all $\zeta\geq q$.
\end{lem}

\begin{proof*}
By construction, the $l^\mathrm{th}$ summand of $M(q)$ in (\ref{C2.7}) has the form 
\linebreak
$\Matrix{b_l}{-a_l}{b_{l-1}}{-a_{l-1}}$ with $\frac{a_{l-1}}{b_{l-1}} < \frac{a_l}{b_l}$.
Since $\frac{a_0}{b_0}=-\infty$ and $\frac{a_L}{b_L}=q$ we find $\frac{a_{l-1}}{b_{l-1}} < q$ for $l = 1,\ldots,L$ and therefore $ b_{l-1}q-a_{l-1} >0$.
Since $\zeta \geq q$ and $c_l=b_{l-1} \geq 0$ the lemma follows immediately.
\end{proof*}

\begin{lem}
\label{C2.9}
Let $A =\Matrix{a}{b}{0}{d} \in \mathrm{Mat}_\ast(2,\Z)$ be such that $a,b \in \N$, $0 \leq b < d$ and $M\big(\frac{b}{d}\big)=\sum_{l=1}^L m_l$.
Then the matrices $m_lA$ contain only nonnegative integer entries.
\end{lem}

\begin{proof*}
By construction, $m_l=\Matrix{b_l}{-a_l}{b_{l-1}}{-a_{l-1}}$ with $\frac{a_{l-1}}{b_{l-1}} < \frac{a_l}{b_l}$ for all $l=1,\ldots,L$ and $\frac{a_0}{b_0}=-\infty$ and $\frac{a_L}{b_L}=\frac{b}{d}$.
Hence the statements $\frac{a_l}{b_l} \leq \frac{b}{d}$ and $a_ld \leq bb_l$ are equivalent for all $l=0,\ldots,L$.
We see that $m_lA=\Matrix{ab_l}{bb_l-da_l}{ab_{l-1}}{bb_{l-1}-da_{l-1}}$ has only nonnegative entries.
\end{proof*}

\begin{lem}
\label{C2.11}
Let $A$ and $m_l$ be as in Lemma~\ref{C2.9}.
Then the entries of the matrix $m_lA=\Matrix{a^\prime}{b^\prime}{c^\prime}{d^\prime}$ satisfy
\[
a^\prime > c^\prime \geq 0 
\qquad \mbox{and} \qquad
d^\prime > b^\prime \geq 0.
\]
\end{lem}

\begin{proof*}
For $q=\frac{b}{d}$ the sequence $\mathrm{LNS}\big(\frac{b}{d}\big)$ in reversed order is minimal according to Lemma~\ref{C2.4}. 
This allows us to use statement (6.3) in \cite{hilgert:1} which is formulated there only for certain upper triangular matrices $A_{[c:d]}$.
The proof of this statement extends however also to the upper triangular matrix $A$.
\end{proof*}

%% file: cusp.tex
\section{Vector valued cusp forms and Hecke operators}
\label{D}
Fix an $n \in \N$ throughout this section.
We will introduce vector valued cusp forms transforming under a representation of the full modular group and show these vector valued cusp forms and scalar valued cusp forms to be equivalent by an explicit bijective map.

\subsection{Hecke operators for scalar valued cusp forms}
\label{D5}
Recall the definition of the Hecke operators on $S(n,s)$ for a fixed $s \in \C$ as introduced by A.~Atkin and J.~Lehner in~\cite{atkin:1}.

\begin{defn}
\label{D5.2}
Denote by $T(p)$ and $U(q)$ for $\gcd(p,n)=1$, $q|n$ and $p,q$ prime the following elements in $\RR_p$ resp.\ $\RR_q$:
\begin{equation}
\label{D5.2a}
\label{D5.2b}
T(p) =
\sum_{ad=p \atop 0 \leq b <d} \Matrix{a}{b}{0}{d} 
\quad \mbox{and} \quad
U(q) =
\sum_{0 \leq b <q} \Matrix{1}{b}{0}{q}.
\end{equation}
The induced maps $S(n,s) \to S(n,s)$ given by $u \mapsto u\big|_0T(p)$ resp.\ $u \mapsto u\big|_0U(q)$ are called the $p^\mathrm{th}$ and $q^\mathrm{th}$ Hecke operator $H_p$ resp.\ $H_q$ on $S(n,s)$
\end{defn}

\noindent
\textit{Remarks.} 
\begin{itemize}
\item
Obviously, the Hecke operators $H_p$ and $H_q$ depend on $n$.
\item
A complete discussion of the Hecke algebra acting on cusp forms for $\Gnull{n}$ can be found in~\cite{miyake:1}. 
\item
For $m\in \N$ the $m^\mathrm{th}$ Hecke operator $H_m$ on $S(1,s)$ is given by
\begin{equation}
\label{D5.1}
H_m u = \sum_{A \in X_m} u\big|_0A,
\end{equation}
see e.g.~\cite{miyake:1}.
\end{itemize}

\subsection{Induced representations}
\label{D1}
Let $G$ be a group and $H$ be a subgroup of $G$ of finite index $\mu=[G:H]$. 
For each representation $\chi:H \to \mathrm{End}(V)$ we consider the induced representation $\chi_H:G \to \mathrm{End}(V_G)$, where
\[
V_G:=\{f:G\to V;\, f(hg)=\chi(h)f(g) \quad \mbox{for all }g\in G, h \in H\} 
\]
and
\[
\big(\chi_H(g)f\big)(g^\prime) = f(g^\prime g)
\qquad \mbox{for all }g,g^\prime \in G.
\]
For $V=\C$ and $\chi$ the trivial representation we call the induced representation $\chi_H$ the \emph{right regular representation}.
In fact, in this case $V_G$ is the space of left $H$-invariant functions on $G$ or, what is the same, functions on $H \backslash G$, and the action is by right translation in the argument. 
One can identify $V_G$ with $V^\mu$ using a set $\{\alpha_1,\ldots,\alpha_\mu\}$ of representatives for $H \backslash G$, i.e.,
\[
H \backslash G = \{H\alpha_1, \ldots, H\alpha_\mu\}.
\]
Then
\[
V_G \to V^\mu 
\quad \mbox{with} \quad
f \mapsto \big(f(\alpha_1),\ldots,f(\alpha_\mu)\big)
\]
is a linear isomorphism which transports $\chi_H$ to the linear $G$-action on $V^\mu$ given by
\[
g \cdot (v_1, \ldots,v_\mu) =  
\big(\chi(\alpha_1 g \alpha_{k_1}^{-1}) v_{k_1}, \ldots,  \chi(\alpha_\mu g \alpha_{k_\mu}^{-1}) v_{k_\mu}\big)
\]
where $k_j\in \{1,\ldots,\mu\}$ is the unique index such that $H\alpha_j g =H\alpha_{k_j}$. 
To see this, one simply calculates 
\[
\big( \chi_H(g)f \big)(\alpha_j) = f(\alpha_j g) = f(\alpha_j g \alpha_{k_j}^{-1} \alpha_{k_j})
= \chi(\alpha_j g \alpha_{k_j}^{-1}) \big(f(\alpha_{k_j})\big).
\]
In the case of the right regular representation the identification $V_G \cong \C^\mu$ gives a matrix realization
\[
\rho(g) = \big( \delta_H(\alpha_i g \alpha_j^{-1}) \big)_{1 \leq i,j \leq \mu}
\]
where $\delta_H(g)=1$ if $g \in H$ and $\delta_H(g)=0$ otherwise.
In particular, the matrix $\rho(g)$ is a permutation matrix.

We take $G=\Gmod$, $H=\Gnull{n}$ and $\alpha_1, \ldots, \alpha_\mu\in \Gmod$ as representatives of the $\Gnull{n}$ orbits in $\Gmod$.
The matrix representation $\rho: \Gmod \to \C^{\mu \times \mu}$ is
\begin{equation}
\label{D1.2}
\rho(g):= 
\Big(  \delta_\Gnull{n} (\alpha_i\,g\, \alpha_j^{-1})  \Big)_{1 \leq i,j \leq \mu}
\qquad \mbox{for all }g \in \Gmod.
\end{equation}
We easily check that $\rho$ satisfies $\rho(g^\prime)\,\rho(g) = \rho(g^\prime g)$ for all $g, g^\prime \in \Gmod$.

\subsection{Vector valued cusp forms}
\label{D2}
For each $u\in S(n,s)$ we construct a vector valued version of $u$ which transforms under the representation $\rho$.
As usual, the index $[\Gmod:\Gnull{n}]$ is denoted by $\mu=\mu_n$.

\begin{defn}
\label{D2.1}
A \emph{vector valued cusp form} $\vec{u}:\HH \to \C^\mu$ for $\Gnull{n}$ with spectral value $s \in \C$ is a vector valued function $\vec{u} = (u_1, \ldots,u_\mu)^\mathrm{tr}$ satisfying
\begin{itemize}
\item $u_j$ is real-analytic for all $j\in \{1,\ldots,\mu\}$,
\item $\vec{u}(g z) = \rho(g)\, \vec{u}(z)$ for all $z \in \HH$ and $g \in \Gmod$,
\item $\Delta u_j = s(1-s)u_j$ for all $j\in \{1,\ldots,\mu\}$ and
\item $u_j(z) = \mathrm{O}\left(\im{z}^C\right)$ as $\im{z} \to \infty$ for all $C \in \R$ and $j \in \{1,\ldots,\mu\}$.
\end{itemize} 
We denote the space of all vector valued cusp forms with spectral parameter $s$ for $\Gnull{n}$ by $S_\mathrm{ind}(n,s)$.
\end{defn}

\textit{Remark.}
The group $\Gmod$ acts on the vector sapce $S_\mathrm{ind}(n,s)$ via
\[
\Gmod \times S_\mathrm{ind}(n,s) \to S_\mathrm{ind}(n,s); \qquad
(g,\vec{u}) \mapsto \rho(g^{-1}) \, \vec{u} \big|_0 g
\]
where $\rho$ is the representation in~(\ref{D1.2}).
This implies in particular that the growth condition on vector valued cusp forms at the cusp $\infty$ gives a growth condition at all cuspidal points $p \in \Gmod \infty = \Q \cup \{\infty\}$.

\medskip

To each $u \in S(n,s)$ we associate the vector valued function $\Pi(u)$ given by
\begin{equation}
\label{D2.2}
\Pi: S(n,s) \to S_\mathrm{ind}(n,s); \quad u \mapsto \Pi(u):= \big( u\big|_0\alpha_1, \ldots, u\big|_0\alpha_n \big)^\mathrm{tr}.
\end{equation}
The function $\Pi(u)$ satisfies all four properties of a vector valued cusp form in Definition~\ref{D2.1}.
Indeed, take an $u \in S(n,s)$ and an index $i \in \{1,\ldots,\mu\}$.
Obviously, $\big[\Pi(u)\big]_i = u\big|_0 \alpha_i$ is real-analytic on $\HH$, it is shown in \S2.2.3 and \S2.2.4 in \cite{bruggeman:13} that
\[
\Delta \big[\Pi(u)\big]_i = \Delta \big( u\big|_0 \alpha_i \big) = (\Delta u)\big|_0 \alpha_i =\big[\Pi(\Delta u)\big]_i,
\]
and the growth condition for $\Pi(u)$ also follows directly from the growth condition for $u$.
To check the transformation property under $\rho$, take a $g \in \Gmod$.
There exists a $g^\prime \in \Gnull{n}$ and an unique $j \in \{1,\ldots,\mu\}$ such that $\alpha_i g \alpha_j^{-1} = g^\prime \in \Gnull{n}$.
Hence
\begin{eqnarray*}
\big[\Pi(u)\big]_i(gz)
&=& 
\big(u\big|_0\alpha_ig \big) (z) = \big(u\big|_0 g^\prime \alpha_j \big) (z)\\
&=&
\sum_{j^\prime=1}^\mu \delta_{\Gnull{n}}\big(\alpha_i g \alpha_{j^\prime}^{-1}\big) \, \big(u\big|_0 g^\prime \alpha_{j^\prime}\big) (z) = \big[\rho(g) \Pi(u|_0g^\prime)\big]_i (z)\\
&=&
\left[ \rho(g)\, \Pi(u) \right]_i (z)
\end{eqnarray*}
The second property follows since $g^\prime \in \Gnull{n}$ and $u$ is $\Gnull{n}$-invariant.
Hence $\Pi(u) \in S_\mathrm{ind}(n,s)$.

On the other hand, consider a vector valued cusp form $\vec{u} \in S_\mathrm{ind}(n,s)$ and take the unique $j\in \{1,\ldots,\mu\}$ with $\Gnull{n}\alpha_j=\Gnull{n}$.
The function $u:=u_j$ is in $S(n,s)$:
The function $u$ satisfies the transformation property $u\big|_0g=u$ for all $g \in \Gnull{n}$ since $\rho(g)_{jj}=\delta_\Gnull{n}(\alpha_j g \alpha_j^{-1}) = 1$ and $u$ is an eigenfunction of $\Delta$ with spectral parameter $s$.
To show that $u$ vanishes in all cusps take a cuspidal point $p \in \Q\cup \{\infty\}$ of $\Gnull{n}$ and $g\in\Gmod$ such that $gp=\infty$.
There exists an index $i \in \{1,\ldots,\mu\}$ and a $\gamma \in \Gnull{n}$ such that $g = \gamma \alpha_i$.
We find that
\[
u_j(\alpha_i z) = \big[\rho(\alpha_i) \vec{u}\big]_j (z) = \big[\vec{u}\big]_i(z)
\]
since $\big[\rho(\alpha_i)\big]_{ji}= \delta_\Gnull{n}\big( \alpha_j\, \alpha_i \, \alpha_i^{-1} \big) = 1$.
Hence
\[
u(gz)= u_j(\gamma \alpha_i z) = u_i(z) = \mathrm{O}\left(\im{z}^C\right)
\quad \mbox{as } \im{z} \to \infty \mbox{ for all } C \in \R.
\]
By the transformation property under the representation $\rho$ one sees that $\vec{u} = \Pi(u)$.

Summarising, we proved the following
\begin{lem}
\label{D2.3}
The spaces $S_\mathrm{ind}(n,s)$ and $S(n,s)$ are isomorphic.
\end{lem}

Moreover, vector valued cusp forms satisfy the following growth condition.
\begin{lem}
\label{D2.4}
For $g \in \Gmod$ a vector valued cusp form $\vec{u} \in S_\mathrm{ind}(n,s)$ satisfies
\begin{equation}
\label{D2.5}
\big[\vec{u}\big]_j(gz) = \mathrm{O}\left(e^{-2\pi \im{z}}\right)
\qquad \mbox{as } \im{z} \to \infty \mbox{ and } \re{z} \mbox{ bounded}
\end{equation}
for all $j \in \{1,\ldots,\mu\}$.
\end{lem}

\begin{proof*}
Take an $u \in S(n,s)$ such that $\Pi(u)=\vec{u}$ and take an index $j \in \{1,\ldots,\mu\}$.
For each cuspidal point $p \in \Q \cup \{i\infty\}$ and $\gamma \in \Gmod$ with $\gamma\infty = p$, Iwaniec has shown in Theorem~3.1 in \cite{iwaniec:1} that
\[
u(\gamma z) = \mathrm{O}\left(e^{-2\pi \im{z}}\right)
\qquad \mbox{as } \im{z} \to \infty \mbox{ and } \re{z} \mbox{ bounded}.
\]
For $j \in \{1,\ldots,\mu\}$ take $p= \alpha_j \infty$ and $\gamma=\alpha_j$.
Hence the above growth estimate implies for the $j^\mathrm{th}$ component of $\vec{u}$ that
\[
\big[\vec{u}\big]_j(z) = \mathrm{O}\left(e^{-2\pi \im{z}}\right)
\qquad \mbox{as } \im{z} \to \infty \mbox{ and } \re{z} \mbox{ bounded}
\]
since $\big[\vec{u}\big]_j = \big[\Pi(u)\big]_j = u\big|_0\alpha_j$.
Repeating the argument for all $j \in \{1,\ldots,\mu\}$ we find that the growth estimate (\ref{D2.5}) holds for $g=I$ the identity element in $\Gmod$.
The stated estimate holds since $\vec{u}(gz) = \rho(g) \, \vec{u} (z)$ for all $g \in \Gmod$ and the index $\mu$ of $\Gnull{n}$ in $\Gmod$ is finite.
\end{proof*}

\subsection{Hecke operators for cusp forms}
\label{D4}
To derive a formula for the Hecke operators acting on $S_\mathrm{ind}(n,s)$ we have to write the vector valued cusp form $\Pi\big( u\big|_0\sum_A A\big)$ in terms of a linear action of a certain matrix sum
 on the vector valued cusp form $\Pi(u)$.

For this recall the Hecke operator $H_m$ in \S\ref{D5}.
For prime $p,q$ with $\gcd(p,n)=1$ and $q|n$ the $p^\mathrm{th}$ (resp.\ $q^\mathrm{th}$) Hecke operator $H_p$ (resp.\ $H_q$) is given by the action of $T(p)$ (resp.\ $U(q)$) on the space of cusp forms which we write as
\[
S(n,s) \to S(n,s); \quad
u \mapsto u \big|_0 \sum_A A = \sum_{A \in \mathbf{A}} u \big|_0 A
\]
with $\sum_A A = T(p)$ and $\mathbf{A} =X_p$ (resp.\ $\sum_A A = U(q)$ and $\mathbf{A}=X_q \smallsetminus \Matrix{q}{0}{0}{1}$).

Consider the $j^\mathrm{th}$ component of the vector valued cusp form $\Pi\big( u\big|_0\sum_A A\big)$.
We would like to write this component as
\begin{equation}
\label{D4.10.a}
\left( u\big|_0\sum_A A\right) \big|_0\alpha_j
=
u\big|_0 \sum_A \left(A \alpha_j \right)
=
u\big|_0 \sum_A \alpha_{\phi_A(j)} \sigma_{\alpha_j}(A)
\end{equation}
for certain indices $\phi_A(j) \in \{1,\ldots,\mu_n\}$ and certain upper triangular matrices $\sigma_{\alpha_j}(A)$.
This will allow us to use the in \S\ref{E} introduced integral transform to determine the form of the Hecke operators on period functions.

The following lemmas show that relation (\ref{D4.10.a}) makes really sence.

\begin{lem}
\label{D4.6}
For each $g \in \Gmod$ there exists a unique bijective map $\sigma_g:X_m \to X_m$ with
$A \,g\, \left( \sigma_g(A)\right)^{-1} \in \Gmod$ for all $A \in X_m$.
\end{lem}

\textit{Remark.}
The inverse of $\sigma_g$ is given by $\sigma_g^{-1} = \sigma_{g^{-1}}$.
Indeed this follows from
\[
\sigma_g^{-1}(A)gA^{-1} \in \Gmod 
\,\, \iff \,\,
\left(\sigma_g^{-1}(A)gA^{-1}\right)^{-1} = Ag^{-1}\left(\sigma_g^{-1}(A)\right)^{-1} \in \Gmod.
\]

For Lemma~\ref{D4.6} we need the following technical result:

\begin{lem}
\label{D4.2}
For any $m \in \N$, any $A \in X_m$ and arbitrary $g \in \Gmod$ there exist unique matrices $A^\prime, A^{\prime\prime} \in X_m$ with
\begin{equation}
\label{D4.2.a}
A \,g\, \left( A^\prime \right)^{-1} \in \Gmod
\quad \mbox{and} \quad
A^{\prime\prime} \,g\, A^{-1} \in \Gmod.
\end{equation}
\end{lem}

\begin{proof*}
The second statement in (\ref{D4.2.a}) follows from the first one by taking the inverse of $A \,g^{-1}\, \left( A^\prime \right)^{-1}$.
We show the first statement in (\ref{D4.2.a}) in two steps.

For $A=\Matrix{a}{b}{0}{d}\in X_m$ with $\gcd(a,b,d)=1$ and $g \in \Gmod$ put $\tilde{A}= \Matrix{a}{\tilde{b}}{0}{d}$ with $\tilde{b}:=b+\big(\gcd(a,b)-1\big)d$.
Then $\gcd(a,\tilde{b})=1$ and Lemma~5.11 in \cite{hilgert:1} applies to $\tilde{A}$:
for given $\tilde{A}$ and $g$ there exists a unique $A^\prime \in X_m$ such that
\[
\tilde{A} \,g\, \left( A^\prime \right)^{-1} = \Matrix{1}{\gcd(a,b)-1}{0}{1} A g\, \left( A^\prime \right)^{-1}\in \Gmod.
\]
Hence (\ref{D4.2.a}) holds for $A$ and $g$.

We extend the result to each $A = \Matrix{a}{b}{0}{d} \in X_m$ and $g \in \Gmod$:
Write $A=\Matrix{a}{b}{0}{d}$ and put $l=\gcd(a,b)$. 
Then the matrix $\Matrix{\frac{a}{l}}{\frac{b}{l}}{0}{d}$ satisfies the assumption in the previous step.
We find a matrix $\Matrix{a^\prime}{b^\prime}{0}{d^\prime} \in X_\frac{m}{l}$ such that
\[
g^\prime:= \Matrix{\frac{a}{l}}{\frac{b}{l}}{0}{d} g \Matrix{a^\prime}{b^\prime}{0}{d^\prime}^{-1}
\in \Gmod.
\] 
Similarly we find a matrix $\Matrix{a^{\prime\prime}}{b^{\prime\prime}}{0}{d^{\prime\prime}} \in X_l$ such that
\[
\Matrix{l}{0}{0}{1} g^\prime \Matrix{a^{\prime\prime}}{b^{\prime\prime}}{0}{d^{\prime\prime}}^{-1} \in \Gmod.
\]
Hence 
\[
\Matrix{a}{b}{0}{d} g \left[\Matrix{a^{\prime\prime}}{b^{\prime\prime}}{0}{d^{\prime\prime}} \Matrix{a^\prime}{b^\prime}{0}{d^\prime} \right]^{-1} \in \Gmod.
\] 
It is shown in the proof of Lemma~5.11 in \cite{hilgert:1} that the entry $a^\prime b^{\prime\prime}+b^\prime d^{\prime\prime}$ of the triangular matrix 
\[
\Matrix{a^\prime}{b^\prime}{0}{d^\prime} \Matrix{a^{\prime\prime}}{b^{\prime\prime}}{0}{d^{\prime\prime}} =
\Matrix{a^\prime a^{\prime\prime}}{a^\prime b^{\prime\prime}+b^\prime d^{\prime\prime}}{0}{d^\prime d^{\prime\prime}} 
\]
is unique modulo $d^\prime d^{\prime\prime}$.
Hence the matrix $\Matrix{a^\prime a^{\prime\prime}}{\hat{b}}{0}{d^\prime d^{\prime\prime}}$
with $0 \leq \hat{b}< d^\prime d^{\prime\prime}$, $\hat{b} \equiv a^\prime b^{\prime\prime}+b^\prime d^{\prime\prime}$ modulo $d^\prime d^{\prime\prime}$ fulfills the requirements stated in the lemma.
\end{proof*}

\begin{proof*}[ of Lemma~\ref{D4.6}]
For $g \in \Gmod$ define the map $\sigma_g:X_m \to X_m$ as $\sigma_g(A) := A^\prime$, where $A^\prime$ is given in (\ref{D4.2.a}) for $A$ and $g$.
Similarly define the map $\sigma_g^\prime:X_m \to X_m$ as $\sigma^\prime_g(A) := A^{\prime\prime}$.
Then the map $\sigma_g^\prime$ is just the inverse of $\sigma_g$:
For $A\in X_m$ put $\tilde{A}:=\sigma_g^\prime \big(\sigma_g(A)\big)$. 
Then by construction $\tilde{A} \in X_m$.
Hence we have two matrices $A$ and $\tilde{A}$ in $X_m$ satisfying
\[
A g \left(\sigma_g(A)\right)^{-1}
\quad \mbox{and} \quad
\tilde{A} g \left(\sigma_g(A)\right)^{-1} \in \Gmod.
\]
Uniqueness in Lemma~\ref{D4.2} shows that $A=\tilde{A}$, i.e.\ 
\[
\sigma_g^\prime \big(\sigma_g(A)\big) =A.
\]
Exchanging the roles of $\sigma_g$ and $\sigma_g^\prime$ above shows also that
\[
\sigma_g \big(\sigma_g^\prime(A)\big) =A.
\]
Hence the map $\sigma_g:X_m \to X_m$ is bijective with inverse $\sigma_g^\prime$.
\end{proof*}

For positive $n,m\in\Z$ the inclusion $\Gnull{mn} \subset \Gnull{n}$ induces a projection map $\Gnull{mn}\backslash \Gmod \to \Gnull{n}\backslash \Gmod$: if $\alpha_1,\ldots,\alpha_{\mu_n}$ are representatives of $\Gnull{n} \backslash \Gmod$ and $\beta_1,\ldots,\beta_{\mu_{mn}}$ are representatives of $\Gnull{mn} \backslash \Gmod$ we can write this map as a map on the indices of the representatives:
\begin{equation}
\label{D4.11.a}
\chi_{mn,n}:\,\{1,\ldots,\mu_{mn}\} \to \{1,\ldots,\mu_n\}; \quad
\end{equation}
with
\begin{equation}
\label{D4.11.b}
\Gnull{mn} \beta_i \subset \Gnull{n} \alpha_{\chi_{mn,n}(i)} 
\qquad \mbox{for all } i\in\{1,\ldots,\mu_{mn}\}.
\end{equation}
holds.

\begin{defn}
\label{D4.8}
Let $\alpha_1, \ldots, \alpha_{\mu_n}$ be representatives of the right cosets of \linebreak $\Gnull{n}$ in $\Gmod$. 
For $A \in X_m$ we define the map
\begin{equation}
\label{D4.8.b}
\phi_A=\phi_{A,n}: \quad \{1,\ldots,\mu_{n}\} \to \{1,\ldots,\mu_n\}; \quad
  i \mapsto \phi_A(i)
\end{equation}
such that 
\begin{equation}
\label{D4.8.c}
A \alpha_i \in \Gnull{n}\alpha_{\phi_A(i)} \,\sigma_{\alpha_i}(A).
\end{equation}
\end{defn}

\textit{Remarks.}
\begin{itemize}
\item
The map $\phi_A$ in Definition~\ref{D4.8} depends on $m$ and the indices $\mu_n =[\Gmod:\Gnull{n}]$ and so on $m$ and $n$.
Usually we write $\phi_{A,n}=\phi_A$ omitting the index $n$ since $n$ is fixed in the entire discussion.
\item
Lemma~\ref{D4.6} implies that the map $\phi_A$ is well defined through relation~(\ref{D4.8.c}).
\end{itemize}

\smallskip

Now we can define Hecke operators for vector valued cusp forms.

\begin{defn}
\label{D4.5}
For $n,m\in \N$, $m$ prime and $s\in \C$ put $\sum_A A := T(m)$ if $m \nmid n$ and put $\sum_A A := U(m)$ if $m \mid n$.
The \emph{$m^\mathrm{th}$ Hecke operator} $H_{n,m}$ on $\vec{u} \in S_\mathrm{ind}(n,s)$ is defined as
\begin{equation}
\label{D4.5.a}
\big(H_{n,m} \vec{u}\big)_{j} \mapsto  \sum_A \, u_{\phi_A(j)} \big|_0 \sigma_{\alpha_j}(A)
\quad \mbox{for } j\in \{1,\ldots,\mu_n\}.
\end{equation}
\end{defn}

\textit{Remark.}
In (\ref{D4.5.a}) we sum over all $A \in X_m$ if $m \nmid n$ and $A \in X_m\setminus\left\{ \Matrix{m}{0}{0}{1} \right\}$ if $m \mid n$.

\medskip

The $m^\mathrm{th}$ Hecke operator $H_{n,m}$ on $S_\mathrm{ind}(n,s)$ corresponds to the $m^\mathrm{th}$ Hecke operator $H_m$ on $S(n,s)$ as the following proposition shows.

\begin{prop}
\label{D4.7}
We have $\Pi(H_m u) = H_{n,m} \Pi(u)$.
\end{prop}

\begin{proof*}
For $\vec{u}=(u_j)_j \in S_\mathrm{ind}(n,s)$ there exists a cusp form $u \in S(n,s)$ with $\vec{u} =\Pi(u)$.
The $m^\mathrm{th}$ Hecke operator $H_m$ acts on $u$ as $u\big|_0\sum_A A$.
Since $\Pi(u) = \left(u\big|_0 \alpha_j\right)_{j\in\{1,\ldots,\mu_n\}}$ we find $\Pi\left(u\big|_0 \sum_A A \right) = \left((u\big|_0 \sum_A A)\big|_0\alpha_j\right)_{j\in\{1,\ldots,\mu_n\}}$.
Since $A\alpha_j \in \Gnull{n} \alpha_{\phi_A(j)} \sigma_{\alpha_j}(A)$ we have
\[
\sum_A u\big|_0 A \alpha_j 
= 
u\big|_0 \sum_A A\alpha_j 
= 
u\big|_0 \sum_A \alpha_{\phi_A(j)} \sigma_{\alpha_j}(A).
\]
Hence $\Pi(H_m u) = H_{n,m} \Pi(u)$.
\end{proof*}

%% file: period.tex
\section{On period functions}
\label{E}
We fix $n\in \N$ and the congruence subgroup $\Gnull{n}$ of $\Gmod$ with index $\mu$ throughout this section.
We recall briefly the definition of period functions in \S\ref{A}.

Recall that a function $f:(0,\infty) \to \C$ is called \emph{holomorphic} if it is locally the restriction of a holomorphic function.

Denote by $\FE^\ast(n,s)$ the space of vector valued functions $\vec{\psi}:(0,\infty) \to \C^\mu$ which are holomorphic in each component and satisfy the \emph{three term equation}
\begin{equation}
\label{E1.1}
\vec{\psi} = \rho(T^{-1})\, \vec{\psi}\big|_{2s}T  + \rho({T^\prime}^{-1})\, \vec{\psi}\big|_{2s} T^\prime
\end{equation} 
where $\rho:\Gmod \to \C^{\mu \times \mu}$ denotes the right regular representation of $\Gnull{n}$ defined in (\ref{D1.2}).
We call such functions $\vec{\psi}$ \emph{period like functions} for $\Gnull{n}$. 
If a period like function $\vec{\psi}=(\psi_i)_i$ satisfies the growth condition
\begin{equation}
\psi_i(z) = \left\{\begin{array}{ll}
\mathrm{O} \left( z^{\max\{0,-2\re{s}\}} \right) \quad& \mbox{as } z \downarrow 0 \mbox{ and}\\
\mathrm{O} \left( z^{\min\{0,-2\re{s}\}} \right) \quad& \mbox{as } z\to \infty.
\end{array} \right.
\end{equation} 
for all $i\in\{1,\ldots,\mu\}$ we call $\vec{\psi}$ a \emph{period function}.
The space of period functions is denoted by $\FE(n,s)$.

\smallskip

\noindent
\textit{Remarks.}
\begin{itemize}
\item
It was shown by J.~Lewis and D.~Zagier in \cite{lewis:2} that for the full modular group the space of period functions $\FE(1,s)$ is isomorphic to $S(1,s)$.
A.~Deitmar and J.~Hilgert generalize this result to submodular groups of finite index in~\cite{deitmar:1}, and hence in particular for the congruence subgroup $\Gnull{n}$.
\item
Put $M=\Matrix{0}{1}{1}{0}$. 
The authors of \cite{hilgert:1} study solutions of 
\begin{equation}
\label{E1.3}
\vec{\psi} = \rho(T^{-1})\, \vec{\psi}\big|_{2s} T  \pm \rho(M{T^\prime}^{-1})\, \vec{\psi}\big|_{2s} (MT^\prime)
\end{equation}
which is Equation~(4.6) in \cite{hilgert:1}.
Solutions $\vec{\psi}$ of equation~(\ref{E1.3}) correspond to eigenfunctions of the transfer operator $\LL_s$ with eigenvalue $\pm1$, see \cite{hilgert:1}.
\end{itemize}

\subsection{Some technical computations}
\label{E2}
To define the integral transform mapping vector valued cusp forms for $\Gnull{n}$ to period functions for $\Gnull{n}$ in \S\ref{E3} we have to recall the function $R_\zeta(z)$ and the $1$-form $\eta(\cdot,\cdot)$ used already by Lewis and Zagier in \cite{lewis:2}, respectively the author in \cite{muehlenbruch:12}.

\smallskip

The function $R_\zeta(z)$ is defined for $x+iy=z \in \HH$ and $\zeta \in \C \smallsetminus \{x\}$ as follows
\begin{equation}
\label{E2.1}
R_\zeta(z) = \frac{y}{(x-\zeta)^2+y^2}.
\end{equation} 
For $\zeta \in \R$ we have $R_\zeta(z) = \frac{\im{z}}{|z-\zeta|^2}$.
Using the relations
\begin{eqnarray*}
&&\im{gz}=\frac{\det g}{|cz+d|^2} \im{z}, \quad
\dd (gz) = \frac{\det g}{(cz+d)^2} \dd z \quad \mbox{and}\\
&&g\zeta-gz = \frac{\det g}{(c\zeta +d)(cz+d)} (\zeta-z)
\end{eqnarray*}
valid for $g=\Matrix{a}{b}{c}{d} \in \GL{\R}$, $\zeta \in \R$ and $z\in \HH$, a straightforward calculation shows that $R_\zeta(z)$ satisfies the transformation formula
\begin{equation}
\label{E2.2}
\frac{|\det g|}{|c\zeta+d|^2 } \,R_{g\zeta}(gz) = R_\zeta(z)
\end{equation} 
for all $g \in \GL{\R}$ and real $\zeta$.
Moreover $R_\zeta^s(z)$ is an eigenfunction of the hyperbolic Laplace operator, \cite{lewis:2},
\begin{equation}
\label{E2.3}
\Delta \big(R_\zeta(z)\big)^s = 
  s(1-s)\, \big(R_\zeta(z)\big)^s
\qquad (s \in \C).
\end{equation} 

For two smooth functions $u,v$ on $\HH$ define the $1$-form $\eta(u,v)$ as in \cite{lewis:2}:
\begin{equation}
\label{E2.4}
\eta(u,v) := \big(v \partial_yu-u\partial_yv \big)dx + \big(u \partial_xv-v\partial_xu \big)dy.
\end{equation} 

The following Lemma is shown in \cite{lewis:2}:
\begin{lem}
\label{E2.5}
If $u$ and $v$ are eigenfunctions of $\Delta$ with the same eigenvalue, then the $1$-form $\eta(u,v)$ is closed.
If $ z \mapsto g(z)$ is any holomorphic change of variables, then the $1$-form satisfies $\eta(u \circ g,v \circ g) =\eta(u,v) \circ g$.
\end{lem}

\subsection{The period functions of vector valued cusp forms}
\label{E3}
We identify $\pm\infty$ with the cusp $i\infty$.
The action of $\Gmod$ on $\HH$ extends naturally to $\HHstar:=\HH \cup \Q \cup \{\infty\}$. 

By a \emph{simple path} $L$ connecting points $z_0,z_1 \in \HHstar$ we understand a piecewise smooth curve which lies inside $\HH$ except possibly for the initial and end point $z_0,z_1$ and is analytic in all points $\HHstar \smallsetminus \HH$ in the sense of \cite{lang:2} on page~58.
Two simple paths $L_{z_0,z_1}$ and $L^\prime_{z_0,z_1}$ are always homotopic, see \cite{lang:2}.
A \emph{path} $L$ connecting points $z_0,z_1 \in \HHstar$ is given by the union of finitely many simple paths $L_n$, $n=1,\ldots,N$ connecting the points $z_{0,n},z_{1,n} \in \HHstar$ such that $z_{0,1}=z_0$, $z_{1,n}=z_{0,n+1}$ and $z_{1,N}=z_1$.
We say that a path $L$ \emph{lies in the first quadrant} resp.\ \emph{in the second quadrant} if $\re{z} \geq 0$ resp.\ $\re{z} \leq 0$ for almost all $z \in L$.
For distinct $z_0,z_1 \in \HHstar \smallsetminus \HH$ the \emph{standard path} $L_{z_0,z_1}$ is the geodesic connecting $z_0$ and $z_1$.
A standard path $L$ is also a simple path.

\begin{defn}
\label{E3.1}
For $\vec{u} \in S_\mathrm{ind}(n,s)$ and $L_{0,\infty}$ the standard path the integral transform $P:S_\mathrm{ind}(n,s) \to \big(C^\omega(0,\infty)\big)^\mu$, with $C^\omega(0,\infty)$ the space of holomorphic functions on $(0,\infty)$, is defined as
\begin{equation}
\label{E3.2}
\big(P\vec{u}\big)_i(\zeta) = \int_{L_{0,\infty}} \eta\big(u_i, R_\zeta^s \big)
\qquad \mbox{for }
\zeta >0, \, i \in \{1,\ldots,\mu\}.
\end{equation}
Formally, we write (\ref{E3.2}) as
\begin{equation}
\label{E3.2a}
\big(P\vec{u}\big)(\zeta) = \int_{L_{0,\infty}} \eta\big(\vec{u}, R_\zeta^s \big).
\end{equation}
\end{defn}

The integrand $\eta\big(u_i, R_\zeta^s \big)(z)$ in (\ref{E3.2}) is of exponentional decay in $0$ and $\infty$ since the cusp form $\vec{u}(z)$ satisfies the growth condition~(\ref{D2.5}).
Hence the integral in (\ref{E3.2}) exists.

\noindent
\textit{Remarks.}
\begin{itemize}
\item
The notation for the integral transform $P$ does not show its dependence on the spectral paramenter $s$.
\item
The convergence of the integral in (\ref{E3.2}) for the full modular group is also shown in \cite{muehlenbruch:12}.
\item
The $1$-forms $\eta(u_i,R_\zeta^s)$ are closed since $u_i$ and $R_\zeta^s$ are eigenfunctions of the hyperbolic Laplacian with the same spectral parameter $s$.
Hence 
\[
\int_{L^\prime} \eta\big(u_i, R_\zeta^s \big)
=
\int_{L_{0,\infty}} \eta\big(u_i, R_\zeta^s \big) = \big[Pu\big]_i(\zeta)
\]
for arbitrary paths $L^\prime$ homotopic to $L_{0,\infty}$.
\item
The function $P\vec{u}(\zeta)$ is a holomorphic function in $\zeta$. 
This is a consequence of the fact that $R_\zeta^s(z)$ is holomorphic in $\zeta$ as can be seen from~(\ref{E2.1}).
\end{itemize}

The function $P\vec{u}$ in (\ref{E3.2a}) has the following transformation property:
\begin{lem}
\label{E3.3}
For $\zeta>0$ and $\gamma = \Matrix{a}{b}{c}{d} \in \Gmod$ with $a,b,c,d \geq 0$ the function $P\vec{u}$ with $\vec{u} \in S_\mathrm{ind}(n,s)$ satisfies
\begin{equation}
\label{E3.4}
(c\zeta +d)^{-2s} \, \rho(\gamma^{-1}) \, P\vec{u}(\gamma \zeta)
= 
\int_{L_{\gamma^{-1}0,\gamma^{-1}\infty}} \eta\big(\vec{u}, R_\zeta^s\big).
\end{equation}
\end{lem}

\textit{Remark.}
Using the slash action we can rewrite (\ref{E3.4}) as
\[
\rho(\gamma^{-1}) \, \left( \big(P\vec{u}\big)\Big|_s\gamma \right)(\zeta)
= 
\int_{\gamma^{-1}L} \eta\big(\vec{u}, R_\zeta^s\big).
\]

\begin{proof*}[ of Lemma~\ref{E3.3}]
Let $\gamma$ satisfy the assumption in the lemma.
By construction $L_{0,\infty}$ and $\gamma^{-1}L=L_{\gamma^{-1}0,\gamma^{-1}\infty}$ are paths in the second quadrant.
Then
\begin{eqnarray*}
(c\zeta+d)^{-2s} \, \rho(\gamma^{-1}) \, P\vec{u}(\gamma \zeta)
&=&
 \int_{L_{0,\infty}}
      \eta\big(\rho(\gamma^{-1})\vec{u}, (c\zeta+d)^{-2s} \,R_{\gamma\zeta}^s \big)\\
&=&
 \int_{L_{0,\infty}}
       \eta\big(\rho(\gamma^{-1})\vec{u}, (c\zeta+d)^{-2s} R_{\gamma\zeta}^s(\gamma\gamma^{-1}\cdot) \big)\\
&=&
 \int_{L_{0,\infty}}
       \eta\big(\vec{u}(\gamma^{-1} \cdot), R_\zeta^s(\gamma^{-1}\cdot) \big)\\
&=&
 \int_{L_{0,\infty}}
       \eta\big(\vec{u}, R_\zeta^s \big) (\gamma^{-1}\cdot)\\
&=&
  \int_{L_{\gamma^{-1}0,\gamma^{-1}\infty}}
     \eta\big( \vec{u}, R_\zeta^s \big)
\end{eqnarray*}
where we made a change of variables and used Lemma~\ref{E2.5}, the transformation property (\ref{E2.2}) of $R_\zeta$ and the property $\vec{u}(\gamma^{-1} z)=\rho(\gamma^{-1})\vec{u}(z)$ of the vector valued cusp form $\vec{u}$.
The substitution is possible since the singularities of the integrand are in $z=c\zeta+d$ and $z=\zeta$. 
Both values are positive under the assumptions of the lemma.
\end{proof*}

\textit{Remark.} 
The assumption in Lemma~\ref{E3.3} that $\gamma$ has no negative entries makes sure that we do not integrate close to the singularities of the integrand. 
We refer to Proposition~42 in \cite{muehlenbruch:12}. 
The corresponding formulation in \cite{lewis:2} is not strong enough.

\begin{lem}
\label{E3.5}
For $\vec{u} \in S_\mathrm{ind}(n,s)$ each component of the function $P\vec{u}$ in (\ref{E3.4}) satisfies the growth conditions
\[
(P\vec{u})_i(\zeta) = \left\{\begin{array}{ll} 
\mathrm{O} \left( \zeta^{\max \big(0,-2\re{s}\big)} \right) \qquad & \mbox{as } \zeta \downarrow 0 \mbox{ and}\\
\mathrm{O} \left( \zeta^{\min \big(0,-2\re{s}\big)} \right) \qquad & \mbox{as } \zeta \to \infty.
\end{array}\right. 
\]
\end{lem}

\begin{proof*}
The lemma follows directly from the proof of Proposition~44 in~\cite{muehlenbruch:12}.
For $n=1$ the growth estimate was also shown in \cite{lewis:2}.
\end{proof*}

\begin{prop}
\label{E3.6}
For $\vec{u} \in S_\mathrm{ind}(n,s)$ the function $P\vec{u}$ is a period function.
\end{prop}

\begin{proof*}
For $\vec{u} \in S_\mathrm{ind}(n,s)$ Lemma~\ref{E3.5} shows that $P\vec{u}$ satisfies the growth conditions for period functions. 
Hence if $P\vec{u}$ satisfies~(\ref{E1.1}) then the Proposition is true.

Consider the path $L_{0,\infty}$.
We have that $L_{0,\infty} = L_{T^{-1}0,T^{-1}\infty} \cup L_{{T^\prime}^{-1}0,{T^\prime}^{-1}\infty}$ and $T$ and $T^\prime $ have only nonnegative matrix entries.
Using Lemma~\ref{E3.3} we find for all $\zeta>0$
\begin{eqnarray*}
P\vec{u}(\zeta) &=& \int_{L_{0,\infty}} \eta\big(\vec{u}, R_\zeta^s\big)\\
&=& 
\int_{L_{T^{-1}0,T^{-1}\infty}}                   \eta\big(\vec{u}, R_\zeta^s\big)
\; + \;
\int_{L_{{T^\prime}^{-1}0,{T^\prime}^{-1}\infty}} \eta\big(\vec{u}, R_\zeta^s\big)\\
&=& \rho(T^{-1}) \, P\vec{u}(T\zeta) \,+ \, \rho({T^\prime}^{-1}) \,(\zeta+1)^{-2s}\, P\vec{u}(T^\prime\zeta).
\end{eqnarray*}
\end{proof*}

\begin{prop}
\label{E3.8}
For $s \in \C \smallsetminus \Z$, $\re{s} > 0$ the operator $P:S_\mathrm{ind}(n,s) \to \FE(n,s)$ is bijective.
\end{prop}

\begin{proof*}
Put $v_0(y) = \frac{1}{\sqrt{y}} u_j(iy)$, $v_1(y) = \frac{\sqrt{y}}{2\pi i} \left(\partial_x u_j\right)(iy)$ with the unique $j \in \{1,\ldots,\mu\}$ such that $\Gnull{n}\alpha_j =\Gnull{n}$.
Let $r\in \C$ satisfy $0 < \re{r} < 2\re{s}$.
Using (\ref{E2.1}) and (\ref{E2.4}) we find that
\begin{eqnarray}
\label{E3.9}
&&
\int_0^\infty \big(P\vec{u}\big)_j(\zeta)\, \zeta^{r-1} \, d\zeta\\
\nonumber
&&\quad=
\int_0^\infty \big(P\vec{u}\big)_j(\zeta)\, \zeta^{r-1} \, d\zeta\\
\nonumber
&&\quad=
i \int_0^\infty\int_0^\infty u_j(iy) \left(\partial_x R_\zeta^s \right)(iy) - R_\zeta^s(iy) \left(\partial_x u_j\right)(iy) \, dy \, \zeta^{r-1} \, d\zeta\\
&&\quad=
\label{E3.11}
i\int_0^\infty u_j(iy) \, \int_0^\infty \left(\partial_x R_\zeta^s \right)(iy)\, \zeta^{r-1}\,d\zeta\,dy \\
\label{E3.10}
&& \qquad
-i\int_0^\infty \left(\partial_x u_j\right)(iy) \, \int_0^\infty R_\zeta^s(iy) \, \zeta^{r-1}\,d\zeta\,dy
\end{eqnarray}
The growth conditions in Lemma~\ref{E3.5} show that the Mellin transform in (\ref{E3.9}) is well defined.
Using the substitution $\zeta \mapsto y\sqrt{t}$ we can compute the inner integrals in (\ref{E3.11}) and (\ref{E3.10}) explicitely:
\begin{eqnarray*}
&&\int_0^\infty \left(\partial_x R_\zeta^s \right)(iy) \, \zeta^{r-1}\,d\zeta
=
\int_0^\infty \, \frac{2s\zeta y^s}{(\zeta^2+y^2)^{s+1}} \zeta^{r-1}\,d\zeta
=
s y^{r-s-1} \, \int_0^\infty \, \frac{t^{\frac{r+1}{2}-1}}{(1+t)^{s+1}}\,d\zeta\\
&&\quad=
s y^{r-s-1} \, B\left({\textstyle\frac{r+1}{2},s+\frac{1-r}{2}} \right)
=
y^{r-s-1} \, \frac{\Gamma \left(\frac{r+1}{2}\right) \Gamma\left(s+\frac{1-r}{2}\right)}{\Gamma(s)}
\end{eqnarray*}
and
\begin{eqnarray*}
\int_0^\infty R_\zeta^s(iy) \, \zeta^{r-1}\,d\zeta
&=&
\frac{1}{2} y^{r-s} \, \int_0^\infty \frac{t^\frac{r}{2}-1}{(1+t)^s} dt\\
&=&
\frac{1}{2} y^{r-s} \, B\left({\textstyle\frac{r}{2},s-\frac{r}{2}} \right)
=
y^{r-s} \frac{\Gamma \left(\frac{r}{2}\right) \Gamma\left(s-\frac{r}{2}\right)}{2\Gamma(s)}
\end{eqnarray*}
where $B(a,b)=\Gamma(a)\Gamma(b)/\Gamma(a+b)=\int_0^\infty \frac{t^{a-1}}{(1+t)^{a+b}} dt$ denotes the beta function (see e.g.\ \cite{AS65} Formula 6.2.1).
Hence the Mellin transform in (\ref{E3.9}) can be written as
\begin{eqnarray}
\label{E3.12}
\int_0^\infty \big(P\vec{u}\big)_j(\zeta)\, \zeta^{r-1} \, d\zeta
&=&\quad
i\frac{\Gamma \left(\frac{r+1}{2}\right) \Gamma\left(s+\frac{1-r}{2}\right)}{\Gamma(s)}
 \,\int_0^\infty y^{r-s-1} \,u_j(iy) \,dy\\
\nonumber
&&-
i\frac{\Gamma \left(\frac{r}{2}\right) \Gamma\left(s-\frac{r}{2}\right)}{2\Gamma(s)}
 \,\int_0^\infty y^{r-s} \left(\partial_x u_j\right)(iy) \,dy\\
\nonumber
&=& \quad
i\frac{\Gamma \left(\frac{r+1}{2}\right) \Gamma\left(s+\frac{1-r}{2}\right)}{\Gamma(s)}
 \,\hat{L}_0\left(u_j,r-s+\frac{1}{2}\right)\\
\nonumber
&& +
\pi\frac{\Gamma \left(\frac{r}{2}\right) \Gamma\left(s-\frac{r}{2}\right)}{\Gamma(s)}
 \,\hat{L}_1\left(u_j,r-s+\frac{1}{2}\right)
\end{eqnarray}
with
\[
\hat{L}_0(u_j,r) = \int_0^\infty u_j(iy) \, y^{r-\frac{3}{2}} dy
\quad\mbox{and}\quad
\hat{L}_1(u_j,r) = \frac{1}{2\pi i} \int_0^\infty \left(\partial_x u_j\right)(iy) \, y^{r-\frac{1}{2}} dy
\]
as defined in Formula~(29) in \cite{deitmar:1}.
The expression at the right hand of (\ref{E3.12}) is (up to a nonzero factor depending on $s$) equal to the expression in Formula~(31) in \cite{deitmar:1}.
Hence $\big[P\vec{u}\big]_j(\zeta)=\Psi_{u_j}(\zeta)$ in \S3 of \cite{deitmar:1}.
Theorem~3.3 in \cite{deitmar:1} implies that the operator $P$ is bijective since $\Pi(u_j)=\vec{u}$.
\end{proof*}

\textit{Remark.}
It is shown in \cite{muehlenbruch:12} that $P\vec{u}$ extends holomorphically to a function on the cut plane $\C \smallsetminus  \R_{\leq 0}$ .

%% file: results.tex
\section{Main result}
\label{F}
In this section we will determine the action of the Hecke operators on the period functions induced from their action on the vector valued cusp forms.
In \S\ref{F1} we discuss the case of the full modular group, which is simpler since we have to deal with scalar valued functions. 
In \S\ref{F2} we discuss the Hecke operators on the period functions for arbitrary $\Gnull{n}$.

\subsection{Hecke operators for period functions for $\Gmod$}
\label{F1}

\begin{lem}
\label{F1.1}
For rational $q \in [0,1)$ put $M(q)=\sum_l m_l$ as in Definition~\ref{C2.6}.
The two paths $L_{q,\infty}$ and $\bigcup_l L_{m_l^{-1} 0, m_l^{-1} \infty}$ have the same initial and end point.  
\end{lem}

\begin{proof*}
Let $M(q)=\bigcup_{l=1}^L m_l$ be defined as in (\ref{C2.7}).
By construction $m_1^{-1}=I$ is the identity matrix and
\[
m_l^{-1} \, 0      = \Matrix{-a_{l-1}}{a_l}{-b_{l-1}}{b_l} \, 0      = \frac{a_l}{b_l}
\quad \mbox{and} \quad
m_l^{-1} \, \infty = \Matrix{-a_{l-1}}{a_l}{-b_{l-1}}{b_l} \, \infty = \frac{a_{l-1}}{b_{l-1}} 
\]
for all $l=1,\ldots,L$. 
In particular $\frac{a_L}{b_L} = q$ and $\frac{a_0}{b_0} = -\infty$.
Hence, we find
\begin{eqnarray}
\label{F1.1.1}
m_L^{-1} \, 0 &=& q,\\
\nonumber
m_l^{-1} \, \infty &=& m_{l-1}^{-1} \, 0 \mbox{ for all } l=1,\ldots,L \quad \mbox{and} \\
\nonumber
m_1^{-1} \, \infty &=& \infty.
\end{eqnarray}
Moreover $\bigcup_{l=1}^L L_{m_l^{-1} 0, m_l^{-1} \infty}$ is a path as union of simple paths and has the same initial and end point as $L_{q,\infty}$.
\end{proof*}

\begin{figure}
\begin{center}
\includegraphics[scale=.85]{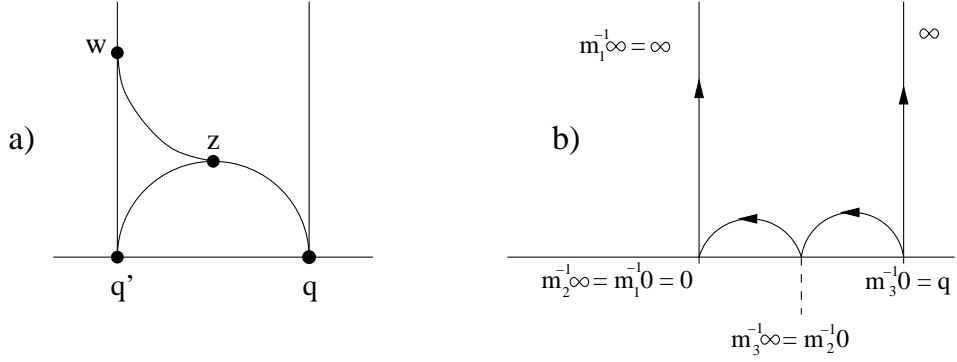}
\caption[Transformation of paths in $\HHstar$]
    {Diagram (a) illustrates the paths used in the proof of Lemma~\ref{F1.7}.
    Diagram (b) illustrates the paths $L_{q,\infty}$ and $\bigcup_l L_{m_l^{-1} 0, m_l^{-1} \infty}$ for e.g.\ $q=\frac{2}{3}$.}
\label{F1.fig1}
\end{center}
\end{figure}

Figure~\ref{F1.fig1}b gives an illustration of the paths $\bigcup_{l=1}^3  L_{m_l^{-1} 0, m_l^{-1} \infty}$ and $L_{q,\infty}$.

\begin{lem}
\label{F1.2}
Let $q$, $L_{q,\infty}$ and $M(q)$ be as in Lemma~\ref{F1.1}.
Both paths $L_{q,\infty}$ and $\bigcup_{l=1}^L L_{m_l^{-1}\,0,m_l^{-1} \,\infty}$ are in the first quadrant.
\end{lem}

\begin{proof*}
Obviously $L_{0,\infty} = \{it; t \geq 0\}$ and $L_{q,\infty} = \{q +it; t \geq 0\}$ ly in the first quadrant since $q \geq 0$. 

The construction of $M(q)$ gives  $m_l^{-1}\,0 = m_{l+1}^{-1} \,\infty \geq 0$, $m_{1}^{-1} \,\infty = \infty$ and $m_L^{-1}\,0 = q \geq 0$.
Hence the paths $L_{m_l^{-1} \,0, m_l^{-1} \,\infty}$ ly in the first quadrant.
The lemma follows by composition of these paths.
\end{proof*}

\begin{lem}
\label{F1.7}
For $\theta$ a closed $1$-form on $\HH$ such that $\int_L \theta$ exists for all simple paths $L$ in $\HHstar$ we have
\[
\int_{L_{q,\infty}} \theta = \int_{L_{q,q^\prime}} \theta + \int_{L_{q^\prime,\infty}} \theta
\]
for all $q,q^\prime \in \Q$.
\end{lem}

\begin{proof*}
Choose points $z$ and $w \in \HH$ on the simple paths $L_{q,q^\prime}$ respectively $L_{q^\prime,\infty}$ as illustrated in Figure~\ref{F1.fig1}a.
We consider the simple paths $L^1 = L_{q^\prime,z} \cup L_{z,w} \cup L_{w,q^\prime}$ with simultanious initial and end point $q^\prime$, $L^2$ with initial point $q$ and end point $\infty$ through the points $z$ and $w$ and $L^3$ with initial point $z$ and end point $\infty$ through the point $w$.

We find that
\[
\int_{L_{q,\infty}} \theta = \int_{L^2} \theta
=
\int_{L_{q,z} \cup L^3} \theta = \int_{L_{q,z}} \theta + \int_{L^3} \theta
\]
since $\theta$ is closed and composition $L^2= L_{q,z} \cup L^3$.
Moreover $\int_{L^1} \theta$ vanishes.
Hence we have
\begin{eqnarray*}
\int_{L_{q,\infty}} \theta 
&=&
\int_{L_{q,z}} \theta + \int_{L^1} \theta + \int_{L^3} \theta\\
&=&
\int_{L_{q,z}} \theta + \int_{L_{z,w}} \theta + \int_{L_{w,q^\prime}} \theta + \int_{L_{q^\prime,z}} \theta + \int_{L^3} \theta\\
&=&
\int_{L_{q,q^\prime}} \theta + \int_{L_{q^\prime,\infty}} \theta.
\end{eqnarray*}
\end{proof*}

We will show in this section how the Hecke operators $H_m$ induce operators $\tilde{H}_m$ on period functions.

\begin{lem}
\label{F1.4}
Let $u \in S(1,s)$ be a cusp form and $Pu$ the period function as in Definition~\ref{E3.1}.
For $A =\Matrix{\ast}{\ast}{0}{d} \in X_m$ and $M(A0) = \sum_{l=1}^L m_l \in \RR_1$ we have
\begin{equation}
\label{F1.4.a}
m^s d^{-2s} \, \int_{L_{A0,A\infty}} \eta\big(u, R_{A\zeta}^s\big) 
 = 
\sum_{l=1}^L \big(Pu \big|_s m_lA \big) (\zeta)= \big(Pu \big|_s M(A0)A\big)(\zeta)
\end{equation}
for all $\zeta >0$.
\end{lem}

\begin{proof*}
For $l \in \{1,\ldots,L\}$ we find
\[
\int_{L_{m_l^{-1}0,m_l^{-1}\infty}} \eta\big(u, R_{A\zeta}^s\big)
 =
\int_{L_{0,\infty}} \eta\big(u\big|_0 m_l^{-1}, R_{A\zeta}^s\big|_0 m_l^{-1}\big)
\]
using Lemma~\ref{E2.5}.
Since $m_l^{-1} \in \Gmod$ and $u$ is a cusp form one has $u\big|_0 m_l^{-1}=u$.
Using~(\ref{E2.2}) we find for $m_l=\Matrix{\ast}{\ast}{c_l}{d_l}$
\[
R_{A\zeta}^s\big|_0 m_l^{-1} 
= 
(c_lA\zeta+d_l)^{-2s}\,R_{m_lA\zeta}^s\big|_0 m_l m_l^{-1}
=
(c_lA\zeta+d_l)^{-2s}\,R_{m_lA\zeta}^s.
\]
The expression $(c_lA\zeta+d_l)^{-2s}$ is well defined since $c_lA\zeta+d_l>0$ (see Lemma~\ref{C2.8}).

Recall that $A=\Matrix{\ast}{\ast}{0}{d}$ and $\det A =m$. 
We have
\begin{eqnarray*}
&& m^s d^{-2s} \, \int_{L_{m_l^{-1}0,m_l^{-1}\infty}} \eta\big(u, R_{A\zeta}^s\big)\\ 
&& \quad = \,
m^s d^{-2s}(c_lA\zeta+d_l)^{-2s}\, \int_{L_{0,\infty}} \eta\big(u, R_{m_lA\zeta}^s\big)\\
&& \quad= \,
\big(Pu\big|_s  m_l A \big) (\zeta).
\end{eqnarray*}
Since Lemma~\ref{C2.9} shows $m_lA \in \mathrm{Mat}_m^+(2,\Z)$, condition (b) in \S\ref{A} holds and the slash action above is well defined.

On the other hand consider the term $\int_{L_{m_L^{-1}0,\infty}} \eta\big(u, R_{A\zeta}^s\big)$.
Lemma~\ref{F1.7} implies that
\[
\int_{L_{m_L^{-1}0,\infty}} \eta\big(u, R_{A\zeta}^s\big)
=
\int_{L_{m_L^{-1}0,m_L^{-1}\infty}} \eta\big(u, R_{A\zeta}^s\big)
+ \int_{L_{m_{L-1}^{-1}0,\infty}} \eta\big(u, R_{A\zeta}^s\big)
\]
since the matrices $m_l$ satisfy (\ref{F1.1.1}).
Iterating the argument gives
\[
\int_{L_{m_L^{-1}0,\infty}} \eta\big(u, R_{A\zeta}^s\big)
=
\sum_{l=1}^L \int_{L_{m_l^{-1}0,m_l^{-1}\infty}} \eta\big(u, R_{A\zeta}^s\big).
\]
Since $L_{A0,A\infty} = L_{m_L^{-1}0,\infty}$ we find
\[
m^s d^{-2s} \, \int_{L_{A0,A\infty}} \eta\big(u, R_{A\zeta}^s\big) 
 =
\sum_{l=1}^L \big(Pu\big|_s  m_lA\big) (\zeta)
\]
for $\zeta >0$.
\end{proof*}

\smallskip

Now, we can define linear operators $\tilde{H}_m$ on period functions.

\begin{defn}
\label{F1.5}
For $m \in \N$ define $\tilde{H}(m) \in \RR_m$ as 
\begin{equation}
\label{F1.5.a}
\tilde{H}(m)=\sum_{d|m \atop 0 \leq b <d} M\big({\textstyle \frac{b}{d}}\big) \Matrix{\frac{m}{d}}{b}{0}{d}.
\end{equation}
For $s \in \C$ the formal sum $\tilde{H}(m)$ induces an operator $\tilde{H}_m$ on the space of holomorphic functions on $(0,\infty)$ through
\begin{equation}
\label{F1.5.b}
\tilde{H}_m f = f\big|_s \tilde{H}(m) 
\qquad (f\in C^\omega(0,\infty)).
\end{equation} 
\end{defn}

\textit{Remark.}
The matrices in the formal sum $\tilde{H}(m)$ in (\ref{F1.5.a}) have only nonnegative integer entries as shown in Lemma~\ref{C2.9}.
Hence $f\big|_s \tilde{H}(m)(\zeta)$ in (\ref{F1.5.b}) is well defined for all $\zeta>0$.

\smallskip

Lemma~\ref{C2.11} implies that the set
\[
\{m_lA; \, A \in X_m, \, M(A0)=\sum_{l=1}^Lm_l\},
\]
which contains all matrices appearing in the formal sum $\tilde{H}(m)$, is a subset of
\[
S_m=\left\{ \Matrix{a}{b}{c}{d}; \, a>c \geq 0, d>b \geq 0 \right\} \subset \mathrm{Mat}_n^+(2,\Z).
\]
It is shown in \cite{hilgert:1} that both sets are indeed equal. (The authors in \cite{hilgert:1} assume that $\gcd(a,b,c,d)=1$ but this restriction is not necessary.)
Hence $\tilde{H}(m) \in \RR_m^+$ is given by $\tilde{H}(m) = \sum_{B\in S_m} B$.

\begin{lem}
\label{F1.6}
Let $Pu$ be the period function of the cusp form $u \in S(1,s)$.
For any $m \in \N$ the operator $\tilde{H}_m$ satisfies
\begin{equation}
\label{F1.6.a}
\Big(\tilde{H}_m (Pu) \Big)(\zeta)
=
\int_0^{i\infty} \eta\big(H_m u, R_\zeta^s \big)
\qquad \mbox{for } \zeta>0.
\end{equation}
\end{lem}

\begin{proof*}
Let $Pu$ be the period function of $u \in S(1,s)$. 
Let $H_m$ be the Hecke operator in~(\ref{D5.1}).
For $A=\Matrix{a}{b}{0}{d} \in X_m$ we have using Lemma~\ref{E2.5} and Equation~(\ref{E2.2})
\begin{eqnarray*}
\int_{L_{0,\infty}} \eta\big( u \big|_0 A, R_\zeta^s\big)
&=&
\int_{L_{A0,A\infty}} \eta\left(u, (R_\zeta^s\big|_0 A^{-1}) \right) \\ 
&=& \left(\frac{\det A}{d^2}\right)^s \,
\int_{L_{A0,A\infty}} \eta\big( u, R_{A\zeta}^s\big).
\end{eqnarray*}
Applying Lemma~\ref{F1.4} we find
\[
\int_{L_{0,\infty}} \eta\big( u \big|_0 A, R_\zeta^s\big)
=
\big(Pu \big|_s M(A0)A\big)(\zeta).
\]
Next we can compute $\int_{L_{0,\infty}} \eta\big( H_m u, R_\zeta^s\big)$ for $H_m u$ in~(\ref{D5.1}) and find
\[
\int_{L_{0,\infty}} \eta\big( H_m u, R_\zeta^s\big)
=
\sum_{A \in X_m}
\Big( Pu\big|_s M\left(A0\right) A \Big)(\zeta)
=
Pu\big|_s \tilde{H}(m) (\zeta)
\]
since $\tilde{H}(m) = \sum_{A\in X_m} M(A0)A$.
\end{proof*}

\begin{prop}
\label{F1.8}
For $u \in S(1,s)$ the period function $Pu \in \FE(1,s)$ satisfies the identity
\[
\big(Pu\big)\big|_s\tilde{H}(m) = P\big(u\big|_0 H(m)\big).
\]
\end{prop}

\begin{proof*}
This follows immediately from Lemma~\ref{F1.6}.
\end{proof*}

\textit{Remark.}
Proposition~\ref{F1.8} shows how the Hecke operators $H_m$ on cusp forms for $\Gmod$ induce Hecke operators on period functions for this group.
These Hecke operators are the same as the operators in \cite{hilgert:1}.
The authors of \cite{hilgert:1} derived the operators using only period functions respectively transfer operators for the groups $\Gnull{m}$.
Another derivation of $\tilde{H}_m$ for $\Gmod$ is given also in \cite{muehlenbruch:13} using a criterion found by Choie and Zagier in \cite{choie:1}.
A similar representation of the Hecke operators has been given by L. Merel in \cite{merel:1}.

\subsection{Hecke operators for period functions for $\Gnull{n}$}
\label{F2}
In this section we extend the above derivation of the Hecke operators for period functions for $\Gmod$ to the congruence subgroups $\Gnull{n}$.

\begin{lem}
\label{F2.1}
For $A =\Matrix{\ast}{\ast}{0}{d} \in X_m$ put $M(A0) = \sum_{l=1}^L m_l \in \RR_1$.
If $P\vec{u}$ is a period function of $\vec{u} \in S_\mathrm{ind}(n,s)$, then
\begin{equation}
\label{F2.1.a}
m^s d^{-2s} \, \int_{L_{A0,A\infty}} \eta\big(\vec{u}, R_{A\zeta}^s\big)
 = 
\sum_{l=1}^L  \rho(m_l^{-1}) \, \Big( P\vec{u} \big|_s m_lA \Big) (\zeta)
\qquad \mbox{for } \zeta>0.
\end{equation}
\end{lem}

\begin{proof*}
Using the transformation property of $\eta$ in Lemma~\ref{E2.5} we get for $l\in \{1,\ldots,L\}$
\[
\int_{L_{m_l^{-1}0,m_l^{-1}\infty}} \eta\big(\vec{u}, R_{A\zeta}^s\big) 
 =
\int_{L_{0,\infty}} \eta\big(\vec{u} \big|_0 m_l^{-1}, R_{A\zeta}^s (m_l^{-1} \cdot)\big).
\]
An argument as in the proof of Lemma~\ref{F1.4} gives
\[
R_{A\zeta}^s\big|_0 m_l^{-1} 
= 
(c_lA\zeta+d_l)^{-2s}\,R_{m_l A\zeta}^s (m_l m_l^{-1} \cdot)
=
(c_lA\zeta+d_l)^{-2s}\,R_{m_l A\zeta}^s
\]
if $m_l=\Matrix{\ast}{\ast}{c_l}{d_l}$. 
Since $c_l A\zeta+d_l>0$ the expression $(c_l A\zeta+d_l)^{-2s}$ is well defined.

Using $\vec{u}\big|_0 m_l^{-1} = \rho(m_l^{-1}) \vec{u}$ for $m_l^{-1} \in \Gmod$ we find
\begin{eqnarray*}
&&
m^s d^{-2s} \, \int_{L_{m_l^{-1}0,m_l^{-1}\infty}} \eta\big( \vec{u}, R_{A\zeta}^s\big)\\
&& \quad = \,
m^s d^{-2s}(c_lA\zeta+d_l)^{-2s}\, \rho(m_l^{-1})\, \int_{L_{0,\infty}} \eta\big(\vec{u}, R_{m_lA\zeta}^s\big)\\
&& \quad = \,
\rho(m_l^{-1})\, \Big(P\vec{u}\big|_s  m_lA \Big) (\zeta).
\end{eqnarray*}

Consider $\int_{L_{A0,A\infty}} \eta\big(\vec{u}, R_{A\zeta}^s\big) (z)$. 
The same argument as in the proof of Lemma~\ref{F1.4} shows
\[
\int_{L_{A0,A\infty}} \eta\big(\vec{u}, R_{A\zeta}^s\big)
 =
\sum_{l=1}^L \int_{L_{m_l^{-1}0,m_l^{-1}\infty}} \eta\big(\vec{u}, R_{A\zeta}^s\big).
\]
But this gives 
\[
m^s d^{-2s} \, \int_{L_{A0,A\infty}} \eta\big(\vec{u}, R_{A\zeta}^s\big)
 =
\sum_{l=1}^L \rho(m_l^{-1})\, \big(P\vec{u}\big|_s  m_lA\big) (\zeta).
\]
\end{proof*}

For $n=1$ Lemma~\ref{F2.1} is just Lemma~\ref{F1.4} since $\rho$ is trivial in this case.

\smallskip

In the following we denote the $i^\mathrm{th}$ component of the vector $\vec{u}$ by $[\vec{u}]_i$.

\begin{lem}
\label{F2.2}
Let $\alpha_1, \ldots, \alpha_\mu$ be representatives of the right coset of $\Gnull{n}$ in $\Gmod$ where $\mu= [\Gmod:\Gnull{n}]$.
Let $P\vec{u}$ be the period function of $\vec{u} \in S_\mathrm{ind}(n,s)$.
For $A\in X_m$, $j\in \{1,\ldots,\mu_{n}\}$ let be $M(\sigma_{\alpha_j}(A)0) = \sum_{l=1}^L m_l \in \RR_1$.
Then the following identity holds for all $j \in \{1,\ldots,\mu\}$ and $\zeta >0$:
\begin{equation}
\label{F2.2.a}
\int_{L_{0,\infty}} \eta\big([\vec{u}]_{\phi_A(j)} \big|_0 \sigma_{\alpha_j}(A), R_\zeta^s \big)
=
\sum_{l=1}^L  \left[ \rho(m_l^{-1}) \, P\vec{u} \right]_{\phi_A(j)} \big|_s m_l \sigma_{\alpha_j}(A)\,(\zeta).
\end{equation}
\end{lem}

\begin{proof*}
Write $\vec{u} = (u_j)_j$.
Using Lemma~\ref{E2.5} and property (\ref{E2.2}) of $R_\zeta$ we find for any $j \in \{1,\ldots,\mu_n\}$ and $\zeta>0$:
\begin{eqnarray*}
&&\int_{L_{0,\infty}} \eta\big(u_{\phi_A(j)} \big|_0 \sigma_{\alpha_j}(A), R_\zeta^s \big)
\, = \,
\int_{L_{\sigma_{\alpha_j}(A)0,\sigma_{\alpha_j}(A)\infty}} \eta\big(u_{\phi_A(j)} , R_\zeta^s \big|_0 \big(\sigma_{\alpha_j}(A)\big)^{-1}\big)\\
&& \quad = \,
m^s d_j^{-2s} \, \int_{L_{\sigma_{\alpha_j}(A)0,\sigma_{\alpha_j}(A)\infty}} \eta\Big(u_{\phi_A(j)} , R_{\sigma_{\alpha_j}(A)\zeta}^s \Big)
\end{eqnarray*}
where $\sigma_{\alpha_j}(A)= \Matrix{\ast}{\ast}{0}{d_j}$ is again in $X_m$.
Take $M(\sigma_{\alpha_j}(A)0)=\sum_{l=1}^L m_l \in \RR_1$ and apply Lemma~\ref{F2.1}. 
We have
\begin{eqnarray*}
\int_{L_{0,\infty}} \eta\big( u_{\phi_A(j)} \big|_0  \sigma_{\alpha_j}(A) , R_\zeta^s \big)
&=&
\Big[ \int_{L_{0,\infty}} \eta\big( \vec{u} \big|_0  \sigma_{\alpha_j}(A) , R_\zeta^s \big) \Big]_{\phi_A(j)}\\
&=&
\left[ 
\sum_{l=1}^L  \rho(m_l^{-1}) \, \Big(P\vec{u} \Big|_s  m_l \sigma_{\alpha_j}(A) \Big)(\zeta)
\right]_{\phi_A(j)} .
\end{eqnarray*}
\end{proof*}

\textit{Remark.}
For $u\in S(n,s)$ and $A \in X_m$ Lemma~\ref{D4.6} implies that
\[
u\big|_0 A\alpha_j = u\big|_0 \alpha_{\phi_A(j)} \sigma_{\alpha_j}(A).
\]
Hence equation~(\ref{F2.2.a}) can be written as
\[
\int_{L_{0,\infty}} \eta\big(\Pi(u\big|_0 A), R_\zeta^s \big)
=
\sum_{l=1}^L  \left[ \rho(m_l^{-1}) \, P\Pi(u) \right]_{\phi_A(j)} \big|_s m_l \sigma_{\alpha_j}(A) \, (\zeta)
\]

\smallskip

Lemma~\ref{F2.2} allows us to derive an explicit formula for the action of the Hecke operators on the period functions of $\Gnull{n}$ induced from the action of these operators on $S_\mathrm{ind}(n,s)$ for this group.

\begin{prop}
\label{F2.3}
Let $\alpha_1, \ldots, \alpha_\mu$ be representatives of the right cosets of \linebreak $\Gnull{n}$ in $\Gmod$.
Let $P\vec{u}$ be the period function of $\vec{u} \in S_\mathrm{ind}(n,s)$.
For $m$ prime take $\mathbf{A} \subset X_m$ such that $T(m)=\sum_{A \in \mathbf{A}} A$ if $\gcd(m,n)=1$ resp.\ $U(m)=\sum_{A \in \mathbf{A}} A$ if $m|n$.
The $m^\mathrm{th}$ Hecke operator $\tilde{H}_{n,m}$ acting on $P\vec{u}$ is given by
\begin{equation}
\label{F2.3.a}
\left[ \tilde{H}_{n,m} \big(P\vec{u}\big)\right]_j=
\sum_{A \in \mathbf{A}} \sum_{l=1}^L  
\left[ \rho(m_l^{-1}) \, P\vec{u} \right]_{\phi_A(j)} \big|_s \big(m_l \sigma_{\alpha_j}(A) \big).
\end{equation}
\end{prop}

\textit{Remark.}
We emphasize that the constant $L$ in (\ref{F2.3.a}) depends on $A$ and that $\mathbf{A} = X_m$ for $\gcd(m,n)=1$ resp.\ $\mathbf{A} = X_m \smallsetminus \left\{\Matrix{m}{0}{0}{1} \right\}$ for $m |n$.

\begin{proof*}[ of Proposition~\ref{F2.3}]
The $m^\mathrm{th}$ Hecke operator $H_{n,m}$ acts on $\vec{u}$ as
\[
\left[ H_{n,m} \vec{u} \right]_j =  \sum_A \, u_{\phi_A(j)} \big|_0 \sigma_{\alpha_j}(A)
\qquad \mbox{for } k \in \{1,\ldots,\mu\}.
\]
Applying Lemma~\ref{F2.2} to both sides then gives formula~(\ref{F2.3.a}).
\end{proof*}

%% file: Hecke.bbl
\begin{thebibliography}{A}

\bibitem[AS65]{AS65}
M.~Abramowitz an d I.~A.~Stegun.
\newblock {\em Handbook of mathematical functions}
\newblock Dover Publications, Inc., 1965.

\bibitem[AL70]{atkin:1}
A.~O.~L. Atkin and J.~Lehner.
\newblock Hecke operators on $\Gamma_0(m)$.
\newblock {\em Math. Ann.} 185 (1970) 134 -- 160.

\bibitem[Br94]{bruggeman:13}
R.~W.~Bruggeman.
\newblock {\em Families of Automorphic forms}, volume 88 in {\em Monographs in Mathematics}
\newblock Birkh\"auser Verlag, 1994.

\bibitem[Ch99]{chang:2}
C.-H. Chang.
\newblock {\em Die {T}ransferoperator-{M}ethode f\"ur {Q}uantenchaos auf den
  {M}odulfl\"achen $\Gamma \backslash \mathbb{H}$}.
\newblock PhD thesis, Mathematisch-Naturwissenschaftliche Fakult\"at der
  Technischen Universit\"at Clausthal, Januar 1999.

\bibitem[CZ93]{choie:1}
{Y.~J.} Choie and {D.} Zagier.
\newblock Rational period functions for $\mathrm{PSL}(2,\mathbb{Z})$.
\newblock In {M.} Knopp and {M.} Sheingorn, editors, {\em A Tribune to {E}mil
  {G}rosswald: Number Theory and related Analysis}, volume 143 of {\em
  Contemporary {M}athematics}, pages 89--108.
\newblock American Mathematical Society, 1993.

\bibitem[DH04]{deitmar:1}\htmladdnormallink{
A.~Deitmar and J.~Hilgert.
\newblock The {L}ewis {C}orrespondence for submodular groups.
\newblock {\em e-arxiv}, 2004.
\newblock {\tt \footnotesize http://arXiv.org/abs/math/0404067}.
}{http://arXiv.org/abs/math/0404067}

\bibitem[HMM05]{hilgert:1}
J.~Hilgert, D.~Mayer and H.~Movasati.
\newblock Transfer operators for {$\Gamma_0(n)$} and the {H}ecke operators for
  period functions of $\mathrm{PSL}(2,\mathbb{Z})$.
\newblock {To appear in \em Math. Proc. Camb. Phil. Soc.} (2005).

\bibitem[Hu94]{hurwitz:1}
A.~Hurwitz.
\newblock Ueber die angen\"aherte {D}arstellung der {Z}ahlen durch rational
  {B}r\"uche.
\newblock {\em Mathematische Annalen} 44 (1894) 417--436.

\bibitem[Iw02]{iwaniec:1}
H.~Iwaniec.
\newblock {\em Spectral Methods of Automorphic Forms}, volume 53 of {\em Graduate Studies in Math.}
\newblock American Mathematical Society, 2002.

\bibitem[La76]{lang:2}
S.~Lang.
\newblock {\em Introduction to modular forms}, volume 222 of {\em Grundl.\ math.\ Wiss.}
\newblock Springer Verlag, Berlin, 1976.

\bibitem[LZ01]{lewis:2}
J.~Lewis and D.~Zagier.
\newblock Period functions for {M}aass wave forms. {I}.
\newblock {\em Ann. of Math.} 153 (2001) 191--258.

\bibitem[Ma01]{martin:1}
F.~Martin.
\newblock {\em P\'eriodes de formes modularires de poids 1},
\newblock Th\`ese doctorate, Paris 7, 2001.

\bibitem[MM]{mayer:8}
D.~Mayer and T.~M\"uhlenbruch.
\newblock {\em From the transfer operator for geodesic flows on modular surfaces to the Hecke operators on period functions of $\Gamma_0(n)$}, 
\newblock In {\em Dynamical {S}ystems: from {A}lgebraic to {T}opoligical {D}ynamics}.
\newblock Proccedings of the ESF-Exploratory Workshop, 5-9 July 2004, Bonn.
\newblock To appear in {\em Contemporary Mathematics},
\newblock American Mathematical Society.

\bibitem[Me94]{merel:1}
L.~Merel.
\newblock Universal Fourier expansions of modular forms.
\newblock In G.~Frey, editor, {\em On Artkins conjecture for odd $2$-dimensional representations}, volume 1585 of {\em Lecture Notes in Math.}
\newblock Springer Verlag, Berlin, 1994.

\bibitem[Mi89]{miyake:1}
T. Miyake.
\newblock {\em Modular {F}orms}.
\newblock Springer-Verlag, 1989.

\bibitem[M\"u03]{muehlenbruch:12}
T. M{\"u}hlenbruch.
\newblock {\em Systems of automorphic forms and period functions}.
\newblock PhD thesis, Utrecht University, September 2003.

\bibitem[Mu04]{muehlenbruch:13}\htmladdnormallink{
T. M\"uhlenbruch.
\newblock Hecke operators on period functions for the full modular group.
\newblock {\em IMRN} 4127-4145 (2004).
}{http://www.hindawi.com/journals/imrn/volume-2004/S1073792804143365.html}
\end{thebibliography}
